\documentclass[a4paper,12pt,fleqn]{article}
\usepackage{amsmath,amssymb,amsthm,graphicx,subfigure,float,caption,epstopdf,tabularx,color, bm,amsfonts,epic}
\usepackage[top=1.25in, bottom=1.25in, left=1.0in, right=1.0in]{geometry}
\usepackage{appendix}
\usepackage{multirow}
\allowdisplaybreaks
\newtheorem{thm}{Theorem}[section]
\newtheorem{lem}{Lemma}[section]

\newtheorem{rmk}{Remark}[section]
\newtheorem*{prf}{Proof}
\numberwithin{equation}{section}

\usepackage{indentfirst}
\graphicspath{{Fig}}
\begin{document}

\title{Structure-preserving algorithms for multi-dimensional fractional
Klein-Gordon-Schr\"{o}dinger equation}

\author{Yayun Fu, \quad Wenjun Cai, \quad Yushun Wang\footnote{Correspondence author. Email: wangyushun@njnu.edu.cn.}\\{\small Jiangsu Key Laboratory for NSLSCS,}\\{\small School of Mathematical Sciences,  Nanjing Normal University,}\\{\small  Nanjing 210023, China}\\}
\date{}
\maketitle

\begin{abstract}
This paper aims to construct structure-preserving numerical schemes for multi-dimensional space fractional Klein-Gordon-Schr\"{o}dinger equation, which are based on the newly developed partitioned averaged vector field methods. First, we derive an equivalent equation, and reformulate the equation as a canonical Hamiltonian system by virtue of the variational derivative of the functional with fractional Laplacian. Then, we develop a semi-discrete conservative scheme via using the Fourier pseudo-spectral method to discrete the equation in space direction. Further applying the partitioned averaged vector field methods on the temporal direction gives a class of fully-discrete schemes that can preserve the mass and energy exactly. Numerical examples are provided to confirm our theoretical analysis results at last.\\[2ex]

\textbf{AMS subject classification:} 35R11, 65M70 \\[2ex]
\textbf{Keywords:} Structure-preserving algorithms; Fractional Klein-Gordon-Schr\"{o}dinger equation; Hamiltonian system; Partitioned averaged vector field methods
\end{abstract}

 %%%%%%%%%%%%%%%%%%%%%%%%%%%%%%%%%%%%%%%%%% begin section 1 introduction %%%%%%%%%%%%%%%%%%%%%%%%%%%%%%%%%%%%%%%%%%%%
\section{Introduction}

Due to the memory and genetic property of fractional operators, fractional differential equations \cite{p3,p4} are more suitable to describe or simulate a variety of scientific and engineering problems with memory and hereditary properties than integer order differential equations, thus they are frequently implemented in many scientific fields, such as physics \cite{p20,p25}, hydrology \cite{p22} and finance \cite{p23}. The fractional Klein-Gordon-Schr\"{o}dinger (KGS) equation \cite{p5} is fractional version of the classical KGS equation, which is introduced by considering a random walk model compounded in Einstein's evolution equation \cite {p6}, and can describe the resonant phenomena
in viscous-elastic materials with long-range dispersion processes and memory effects. %the fractional Klein-Gordon-Schr\"{o}dinger (KGS) equation , which is fractional version of classical KGS equation. %Nowadays,  For instance,
%Recently, the investigation for this equation has attracted researchers' attention. For instance, in Ref \cite{p5},
 In Ref. \cite{p5}, Guo studied the global well-posedness of the fractional KGS system, and pointed out that the fractional order of Laplacian reveals the order of dispersive effect.
As is known to all, it is generally difficult to obtain the explicit forms of the analytical solutions of fractional differential equations. In the last few years, many numerical methods for fractional differential equations have been developed by researchers, such as finite difference methods \cite{p26,p27}, spectral methods \cite{p28,p32} and finite element methods \cite{p30,p31}. For the fractional KGS equation with different boundary conditions, numerical schemes were also constructed, which can be referred to Refs. \cite{p7, p8, p35}.

In the paper, we numerically consider the following fractional KGS equation
\begin{align}\label{FKGS:eq:1.1}
\text{i}\varphi_{t}-\frac{1}{2}(-\Delta)^{\frac{\alpha}{2}}\varphi+ u\varphi=0,\ \ \ \ \ \ \  0<t\leq T,
\end{align}
\begin{align}\label{FKGS:eq:1.2}
u_{tt}+(-\Delta)^{\frac{\beta}{2}}u+u-|\varphi|^2=0,\ \ \  0<t\leq T,
\end{align}
 with the initial conditions
\begin{align}\label{FKGS:eq:1.3}
\varphi(\bm x,0)=\varphi_{0}(\bm x),\ \ u(\bm x,0)=u_{0}( \bm x),\ \ u_{t}(\bm x,0)=\tilde{u}_{0}(\bm x),
\end{align}
where $\text{i}^2=-1$, $1< \alpha, \beta \leq 2$, $\varphi(\bm x,t)$ and $u(\bm x,t)$ are complex and real valued functions of $\bm x\in \Omega \subset \mathbb{ R}^d$ ($d$=1, 2, 3), respectively, and subject to the periodic boundary conditions. $\varphi_{0}(\bm x),~u_{0}(\bm x)$ and $\tilde{u}_{0}(\bm x)$ are given smooth functions. Given $1< s \leq 2$, the fractional Laplacian $(-\Delta)^{\frac{s}{2}}$
is defined as a pseudo-differential operator with the symbol $|\bm \xi|^{s}$ in the Fourier space \cite{p36}
\begin{align}\label{NLS:eq:2.3}
\widehat{(-\Delta)^{\frac{s}{2}}}u(\bm \xi)=|\bm \xi|^{s}\widehat{u}(\bm \xi)
\end{align}
where $\widehat{u}(\bm \xi)=\int_{-\infty}^{\infty}u(\bm x)e^{-\text{i}\bm \xi \bm x}d\bm x$ denotes the Fourier transform of $u(\bm x)$.  If $\alpha=\beta=2$, the equation reduces to the standard KGS equation, which has been studied by scholars \cite{p29,p33,p34}.

As well known, many continuous systems possess some physical quantities that naturally arise from the
physical context, such as energy, momentum and mass.
For system \eqref{FKGS:eq:1.1}-\eqref{FKGS:eq:1.2} with periodic boundary conditions, Guo et al. \cite{p5,p8} derived that the system has fractional mass and energy conservation laws
 \begin{align}\label{FKGS:eq:1.7}
{M}(t)={M}(0),\ \ \ {E}(t)={E}(0),
\end{align}
where the mass has the form
\begin{align}\label{FKGS:eq:1.8}
{M}(t):=\int_{\Omega}|\varphi|^2d\bm x,
 \end{align}
 and the energy is defined as
 \begin{align}\label{FKGS:eq:1.9}
{E}(t):=\int_{\Omega}\Big(u_{t}^2+u^2+((-\Delta)^{\frac{\alpha}{4}}\varphi)^2+((-\Delta)^{\frac{\beta}{4}}u)^2-2u|\varphi|^2\Big)d\bm x.
 \end{align}

Structure-preserving algorithms are methods that can conserve at least one of the intrinsic properties of a given dynamical system, prior researches generally confirmed that the methods are more superior due to their excellent stability and accurate long-time simulation \cite{p9,p10}. Nowadays, structure-preserving numerical schemes for fractional differential equations exert a tremendous fascination on scholars \cite{p11,p12,p13,p15}. The conservations of energy and mass are the crucial properties of the fractional KGS equation \cite{p6}, therefore it is of interest to investigate corresponding conservative schemes. Recently, some scholars studied structure-preserving numerical schemes for the one-dimensional fractional KGS equation. For instance, Wang and Xiao developed a conservative difference scheme \cite {p7} for the equation with Dirichlet boundary condition. In Ref. \cite {p8}, they further proposed Fourier spectral method to solve the equation with periodic boundary condition, and proved that the scheme can conserve the discrete mass and energy. Li \cite{p35} constructed fast conservative finite element schemes to solve the equation. However, to the best of our knowledge, structure-preserving algorithms have not been considered for multi-dimensional fractional KGS equation.
  %and constructed conservative schemes for Hamiltonian ODEs and PDEs based on these methods were  in Ref. \cite{p2}. However, up to now, there is no previous research focusing on the PAVF methods for the fractional Hamiltonian PDEs. The spectral methods play an important role in numerically and theoretically solving differential equations, and its main properties are high accuracy and much more efficient.
%In this paper, with the fractional KGS equation as an example, the PAVF methods are extended to solve fractional differential equations. We construct a class of conservative schemes for the fractional KGS equation based on the PAVF method and the Fourier pseudo-spectral method.

The averaged vector field (AVF) method can exactly preserve the energy for Hamiltonian systems \cite{p16, p17},  which can be applied to construct conservative schemes for the energy conservation fractional partial differential equations. However, the schemes are fully-implicit, thus a nonlinear iteration has to be devised. In addition, since all the unknown variables are averaged simultaneously in the numerical schemes, the scale of number of iterations for the nonlinear scheme may be large.
Recently, to reduce iterative costs at each temporal step and maintain the desired energy-preserving property, Cai et al. \cite{p2} developed partitioned averaged vector field (PAVF) methods that can conserve the energy and mass, and have been used to construct conservative schemes for Hamiltonian ODEs \cite {p2}. Compared with the AVF method, the resulted schemes of the PAVF method are much simpler according to the concrete expressions which reduces the original fully implicit schemes to semi-implicit or linearly implicit schemes and therefore significantly improve the computational efficiency. The spectral methods play an important role in numerically and theoretically solving differential equations, and its main properties are high accuracy, and more efficient due to the fast Fourier transformation technique can be used.

The motivation of this paper is to develop more efficient structure-preserving schemes
for multi-dimensional fractional KGS equation. According to the theories of structure-preserving algorithms,
the derivation of the Hamiltonian formulation for the equation is the key process for constructing numerical schemes. However, there is few research focusing on the Hamiltonian formulation of fractional differential equations. Recently, Wang and Huang \cite{p1} presented the variational derivative of the functional with fractional Laplacian, and reformulated the one-dimensional fractional nonlinear Schr\"{o}dinger equation as a Hamiltonian system. Based on this, we further investigate the properties of
the fractional Laplacian and  derive the Hamiltonian formulation of multi-dimensional fractional KGS equation, and then construct conservative schemes for the equation based on the PAVF methods and the Fourier pseudo-spectral method.

The remainder of this paper is arranged as follows. In Section 2, we reformulate multi-dimensional fractional KGS equation as an infinite-dimensional canonical Hamiltonian system. In Section 3, we first derive a semi-discrete conservative scheme by using the Fourier pseudo-spectral approximation to fractional Laplacian operator. Then, we obtain a class of fully-discrete conservative schemes by utilizing the PAVF methods to discrete the semi-discrete system in time. Numerical examples are presented in Section 4 to demonstrate the theoretical results. Some
conclusions are drawn in Section 5.

%%%%%%%%%%%%%%%%%%%%%%%%%%%%%%%%%%%%%%%%%% begin section 1 introduction %%%%%%%%%%%%%%%%%%%%%%%%%%%%%%%%%%%%%%%%%%%%

\section{Hamiltonian formulation and the PAVF methods}

The Hamiltonian structure is important to analyse the conservative systems and further
to construct numerical schemes for them.  Up to now, only a few
researchers consider the Hamiltonian formulation of one-dimensional fractional differential equations \cite{p1}. In this section, we will investigate the properties of
the fractional Laplacian and derive the Hamiltonian formulation of system \eqref{FKGS:eq:1.1}-\eqref{FKGS:eq:1.2}.
Before embarking on this, we first give some auxiliary lemmas with respect to the fractional Laplacian.

\begin{lem}\rm Given $s >0$, then for any real periodic functions $p, q \in L^2({\mathbb{R}})$, we have
\begin{align}\label{NLS:eq:2.4}
 \big((-\Delta)^{s}p,{q}\big)=\big((-\Delta)^{\frac{s}{2}}p,{(-\Delta)^{\frac{s}{2}}q}\big)=\big(p,{(-\Delta)^{s}q}\big).
\end{align}
\end{lem}
\begin{prf}\rm  First, let us recall a useful property of Fourier transform, namely,
\begin{align}\label{NLS:eq:2.5}
\int_{\Omega}pqd{\bm x}=\int_{-\infty}^{\infty}\hat{p}\hat{q}d{\bm \xi}.
\end{align}
Then, we can deduce
\begin{align}\label{NLS:eq:2.6}
 \big((-\Delta)^{s}p,{q}\big)&= \big(\widehat{(-\Delta)^{s}p},\widehat{q}\big)=\big(|\bm \xi|^{2s}\widehat{p},\widehat{q}\big)=\big(|\bm \xi|^{s}\widehat{p},|\bm \xi|^{s}\widehat{q}\big)\nonumber\\
&=\big(\widehat{(-\Delta)^{\frac{s}{2}}p},\widehat{(-\Delta)^{\frac{s}{2}}q}\big)=\big({(-\Delta)^{\frac{s}{2}}p},{(-\Delta)^{\frac{s}{2}}q}\big),
 \end{align}
and
\begin{align}\label{NLS:eq:2.7}
\big((-\Delta)^{s}p,{q}\big)&= \big(\widehat{(-\Delta)^{s}p},\widehat{q}\big)=\big(|\bm \xi|^{2s}\widehat{p},\widehat{q})=(\widehat{p},|\bm \xi|^{2s}\widehat{q}\big)\nonumber\\
&=\big(\widehat{p},\widehat{(-\Delta)^{s}q}\big)=\big(p,{(-\Delta)^{s}q}\big).
\end{align}
\qed
\end{prf}

\begin{lem}\rm
For a functional $F[g]$ with the following form
\begin{align}\label{NLS:eq:2.8}
F[g]=\int_{\Omega}f\big( g({\bm x}), (-\Delta)^{\frac{\alpha}{4}}g({\bm x})\big)d{\bm x},
\end{align}
where $f$ is a smooth function on the $\Omega$, the variational derivative of $F[\rho]$ is given as follows
\begin{align}\label{NLS:eq:2.9}
\frac{\delta F}{\delta g}=\frac{\partial f}{\partial g}+(-\Delta)^{\frac{\alpha}{4}}\frac{\partial f}{\partial\big((-\Delta)^{\frac{\alpha}{4}}g\big)}.
\end{align}
\end{lem}
\begin{prf}\rm
Let $\phi({\bm x})$ be an arbitrary function with the periodic boundary condition. According to the fact that the fractional Laplacian is linear operator, and the definition of variational derivative, we have
\begin{align}\label{NLS:eq:2.10}
\int_{\Omega}\frac{\delta F}{\delta g}\phi({\bm x})d{\bm x}&=\Big[\frac{d}{d\mu}\int_{\Omega}f\big(g+\mu\phi,(-\Delta)^{\frac{\alpha}{4}}g+\mu(-\Delta)^{\frac{\alpha}{4}}\phi \big)d{\bm x}\Big]_{\mu=0}\nonumber\\
&=\int_{\Omega}\big(\frac{\partial f}{\partial g} \phi+\frac{\partial f}{\partial((-\Delta)^{\frac{\alpha}{4}}g)}(-\Delta)^{\frac{\alpha}{4}}\phi \big)d{\bm x} \nonumber\\
&=\int_{\Omega}\Big(\frac{\partial f}{\partial g} \phi+\big((-\Delta)^{\frac{\alpha}{4}}\frac{\partial f}{\partial((-\Delta)^{\frac{\alpha}{4}}g)}\big)\phi \Big)d{\bm x}
\nonumber\\
&=\int_{\Omega}\Big(\frac{\partial f}{\partial g}+(-\Delta)^{\frac{\alpha}{4}}\frac{\partial f}{\partial((-\Delta)^{\frac{\alpha}{4}}g)} \Big) \phi d{\bm x},
\end{align}
where \eqref{NLS:eq:2.4} was used. Based on the fact that $\phi({\bm x})$ is arbitrary, by using the fundamental lemma of calculus of variations, we can obtain \eqref{NLS:eq:2.9}.
\qed
\end{prf}

\subsection{ Hamiltonian formulation}

By setting $\varphi=q+\text{i}p$, $v=\frac{1}{2}u_{t}$, system \eqref{FKGS:eq:1.1}-\eqref{FKGS:eq:1.2} can be rewritten as a first-order system
\begin{align}\label{FKGS:eq:2.4}
u_{t}=2v,
\end{align}
\begin{align}\label{FKGS:eq:2.5}
v_{t}=\frac{1}{2}\Big(-(-\Delta)^{\frac{\beta}{2}}u-u+(p^2+q^2)\Big),
\end{align}
\begin{align} \label{FKGS:eq:2.6}
p_{t}=-\frac{1}{2}(-\Delta)^{\frac{\alpha}{2}}q+ uq,
\end{align}
\begin{align}
\label{FKGS:eq:2.7}
q_{t}=\frac{1}{2}(-\Delta)^{\frac{\alpha}{2}}p-up,
\end{align}
where $u,v,p,q$ subject to the periodic boundary conditions.

\begin{thm}\rm\label{2SG-lem2.1} The system \eqref{FKGS:eq:2.4}-\eqref{FKGS:eq:2.7} with the periodic boundary conditions has the following global invariants
\begin{align*}
\mathcal{ M}&=\int_{\Omega}(p^2+q^2)d\bm x,
\end{align*}
\begin{align*}
\mathcal{H}&=\frac{1}{4}\int_{\Omega}\Big[ \Big((-\Delta)^{\frac{\alpha}{4}}p\Big)^2+\Big((-\Delta)^{\frac{\alpha}{4}}q\Big)^2+\Big((-\Delta)^{\frac{\beta}{4}}u\Big)^2+u^2+4v^2-2 u(p^2+q^2)\Big]d\bm x.
\end{align*}
\end{thm}
\begin{prf}\rm Computing the inner product of \eqref{FKGS:eq:2.6}-\eqref{FKGS:eq:2.7} with $p$ and $q$, respectively, we obtain the conservation of mass
\begin{align*}
\frac{d}{d t}\mathcal{M}=0.
\end{align*}
Then, by taking the inner products of \eqref{FKGS:eq:2.4}-\eqref{FKGS:eq:2.7} with $v_{t}, u_{t}, q_{t}, -p_t$, respectively, one can deduce that the conservation of energy
\begin{align*}
\frac{d}{d t}\mathcal{H}=0,
\end{align*}
where Lemma 2.1 was used.
We finish the proof.
\qed
\end{prf}

Based on the fractional variational  derivative formula in Lemma 2.2, we obtain the following result straightforwardly.

\begin{thm}\rm\label{2SG-lem2.2} The system \eqref{FKGS:eq:2.4}-\eqref{FKGS:eq:2.7} is an infinite-dimensional canonical Hamiltonian system
\begin{align*}
\left(\begin{array}{c}
		v_t\\
		p_t\\
        u_t\\
		q_t
		\end{array} \right)=J^{-1} \left(\begin{array}{c}
		\delta\mathcal{H}/\delta v\\
		\delta\mathcal{H}/\delta p\\
        \delta\mathcal{H}/\delta u\\
		\delta\mathcal{H}/\delta q
		\end{array} \right),\ \ J=\left(\begin{array}{ccccc}
              {0}& 0&1&0\\
              {0}& 0&0&1 \\
              {-1}& 0&0&0\\
              {0}& -1&0&0
              \end{array}
\right)
\end{align*}
where the energy functional $\mathcal{H}$ is given in Theorem 2.1.

\end{thm}
\subsection{Partitioned averaged vector field methods }
 The paper aims to construct conservative schemes for the fractional KGS equation based on the PAVF methods. Therefore, in this subsection, we briefly introduce them.

Considering the following Hamiltonian system
\begin{align}\label{FKGS:eq:3.1}
\frac{dw}{dt}=f(w)=S_{k}\nabla H(w),\ \ w(0)=w_{0},
\end{align}
where $w\in \mathbb{R}^k$, $S_{k}$ is a $k\times k$ skew symmetric matrix, $k$ is an even number, and the Hamiltonian
$H(w)$ is assumed to be sufficiently differentiable. The so-called second-order AVF method for system \eqref{FKGS:eq:3.1} is defined by
\begin{align}\label{FKGS:eq:3.2}
\frac{w^{m+1}-w^{m}}{\tau}&=\int_{0}^{1}f(\varepsilon w^{m+1}+(1-\varepsilon)w^{m} )d\varepsilon\nonumber\\&=S_{k}\int_{0}^{1}\nabla H(\varepsilon w^{m+1}+(1-\varepsilon)w^{m} )d\varepsilon,
\end{align}
where $\tau$ is the time step.

%The partitioned averaged vector field methods were developed by Cai and his collaborators \cite{p2}, we briefly introduce them as follows

Let $k=2d$, $w=(w_{1},w_{2}, \cdots, w_{d}; w_{d+1}, w_{d+2}, \cdots, w_{k})^{T}=(y,z)^{T}$, then the equivalent system of original problem
\eqref{FKGS:eq:3.1} reads as follows
\begin{align}\label{FKGS:eq:3.3}
\left(\begin{array}{ccccc}
              \dot{y}\\
              \dot{z}
              \end{array}
\right)=S_{2d}\left(\begin{array}{ccccc}
              H_{y}(y,z)\\
              H_{z}(y,z)
              \end{array}
\right), \ \ y, z \in \mathbb{R}^{d},
\end{align}
and the Hamiltonian $H(y,z)$ still retains conserved along any continuous flow.

Obviously, the numerical schemes obtained by the AVF method are fully-implicit, to reduce the computational complexity and preserve the Hamiltonian of system, we developed the PAVF method \cite{p2}
for the Hamiltonian system \eqref{FKGS:eq:3.3}, namely
\begin{align}\label{FKGS:eq:3.4}
\frac{1}{\tau}\left(\begin{array}{ccccc}
              {y^{m+1}-y^{n}}\\
              {z^{m+1}-z^{n}}
              \end{array}
\right)=S_{2d}\left(\begin{array}{ccccc}
              \int^{1}_{0}H_{y}(\varepsilon y^{m+1}+(1-\varepsilon)y^{m},z^{m})d\varepsilon\\
              \int^{1}_{0}H_{z}(y^{m+1},\varepsilon z^{m+1}+(1-\varepsilon)z^{m})d\varepsilon
              \end{array}
\right).
\end{align}
One can easily proved that the PAVF method preserves the Hamiltonian $H(y,z)$ of the system \eqref{FKGS:eq:3.3} exactly.
\begin{rmk} \rm
 The key strategy of the PAVF method is to separate variables in cross-product terms into different groups and apply the AVF method only for one group during each step, which lead to much simper schemes that usually require less computational effort than the conventional AVF method. Particularly, the corresponding schemes can be linearly implicit at least for a part of the entire equation system which dramatically reduces the iteration scale or even avoid the nonlinear iteration completely.  Two methods both achieve second order accuracy and energy-preserving property.
\end{rmk}

We note that the PAVF method \eqref{FKGS:eq:3.4} is only of first-order accuracy \cite{p2}. To increase the accuracy and preserve the Hamiltonian, we combine the PAVF method and its adjoint method, then obtain the PAVF composition (PAVF-C) method and PAVF plus (PAVF-P) method.

Let above PAVF method \eqref{FKGS:eq:3.4} be $\Phi_{\tau}$, and its adjoint method $\Phi^{*}_{\tau}$ is defined as follows
\begin{align}\label{FKGS:eq:3.5}
\frac{1}{\tau}\left(\begin{array}{ccccc}
              {y^{m+1}-y^{m}}\\
              {z^{m+1}-z^{m}}
              \end{array}
\right)=S_{2d}\left(\begin{array}{ccccc}
              \int^{1}_{0}H_{y}(\varepsilon y^{m+1}+(1-\varepsilon)y^{m},z^{m+1})d\varepsilon\\
              \int^{1}_{0}H_{z}(y^{m},\varepsilon z^{m+1}+(1-\varepsilon)z^{m})d\varepsilon
              \end{array}
\right).
\end{align}
Then, the PAVF composition (PAVF-C) method is given by
\begin{align}\label{FKGS:eq:3.6}
\Upsilon_{\tau}:=\Phi^{*}_{\frac{\tau}{2}}\circ \Phi_{\frac{\tau}{2}},
\end{align}
and the PAVF plus (PAVF-P) method is defined by
\begin{align}\label{FKGS:eq:3.7}
\hat{\Upsilon}_{\tau}:=\frac{1}{2}(\Phi^{*}_{\frac{\tau}{2}}+ \Phi_{\frac{\tau}{2}}).
\end{align}

\begin{rmk} \rm
Owning to the great advantage of the PAVF method, even by its composition the PAVF-C method still has much lower cost than the direct AVF method. Although the computational efficiency of the PAVF-P method is comparable to the AVF method, it may possess additional conservative quantities that the AVF method cannot preserve.
\end{rmk}

\section{ Construction of the conservative schemes }

In this section, to present the construction process of numerical schemes clearly, we take the one-dimensional case of system \eqref{FKGS:eq:1.1}-\eqref{FKGS:eq:1.2} as an example, i.e., $d=1$, and develop a class of conservative schemes for the system. Similarly, the construction processes of the schemes can be generalized to the multi-dimensional case.
First, we truncate into a finite computational domain $\Omega=(x_{a}, x_{b})$. The fractional Laplacian $-(-\Delta)^{\frac{s}{2}}$, under the bounded interval $\Omega$ with the periodic boundary condition, is defined by the Fourier series as \cite {p1}
\begin{align}\label{FKGS:eq:1.5}
-(-\Delta)^{\frac{s}{2}}u(x,t)=-\sum\limits_{k\in \mathbb{Z}}|\upsilon_{k}|^{s}\hat{u}_{k}e^{\text{i}\upsilon_{k}(x-x_{a})},
\end{align}
where $\upsilon_{k}=\frac{2k\pi}{x_{b}-x_{a}}$, $\hat{u}_{k}$ is Fourier coefficient and is given by
\begin{align}\label{FKGS:eq:1.6}
\hat{u}_{k}=\frac{1}{x_{b}-x_{a}}\int_{\Omega}u(x,t)e^{-\text{i}\upsilon_{k}(x-x_{a})} dx.
\end{align}

\subsection{Structure-preserving spatial discretization}

We define the time step $\tau:= \frac{T}{M}$ and mesh size $h:=\frac{x_{b}-x_{a}}{N}$ with the integers $M$ and $N$, and denote~$\Omega_{h}=\{x_j |\ x_j=x_a+jh, 0\leq j \leq N-1 \}$,
$\Omega_{\tau}=\{t_m|\ t_m=m\tau, 0\leq m \leq M  \}$. Let $(u_{j}^{m},v^{m}_{j},p_{j}^{m},q^{m}_{j})$ be the numerical approximations to the exact solutions $(u,v,p,q)$ at the grid point $(x_{j},t_{m})$. The corresponding vector forms at any time level are then given by
\begin{align*}
&U=(u_0,u_1, \cdots, u_{N-1})^T,\ \ \ V=(v_0,v_1, \cdots, v_{N-1})^T, \\ &P=(p_0,p_1, \cdots, p_{N-1})^T,\ \ \ Q=(q_0,q_1, \cdots, q_{N-1})^T.
\end{align*}
With these notations, we define the discrete inner product and the associated discrete maximum norm ($l^{\infty}$-norm) as
 \begin{align*}
(P, Q)=h\sum\limits_{j=0}^{N-1} p_{j}{{q}}_{j}, \ \ \ \|P\|=(P,P)^{\frac{1}{2}}, \ \ \|P\|_{{\infty}}=\sup\limits_{0 \leq j \leq N-1}|p_j|,
 \end{align*}
and introduce following operators
\begin{align*}
\delta_{t}P^{m}=\frac{P^{m+1}-P^{m}}{\tau},\ \ \ P^{m+\frac{1}{2}}=\frac{P^{m+1}+P^{m}}{2}.
\end{align*}

Let $x_{j}\in \Omega_{h}$ be the Fourier collocation points. Using the interpolation polynomial $I_{N}{u}(x)$ to approximate ${u}_{N}(x)$ of the function ${u}(x)$, and which is defined by
\begin{align*}
I_{N}{u}(x)={u}_{N}(x)=\sum\limits_{k=-N/2}^{N/2}\hat{u}_{k}e^{\text{i}k\mu(x-x_{a})},
\end{align*}
where $\mu=\frac{2\pi}{x_{b}-x_{a}}$, and the coefficient
\begin{align*}
\hat{u}_{k}=\frac{1}{Nc_{k}}\sum\limits_{j=0}^{N-1}{u}(x_{j})e^{-\text{i}k\mu(x-x_{a})},
\end{align*}
where $c_{k}=1$ for $|k|< N/2$, and $c_{k}=2$ for $k=\pm N/2$. Then
 the fractional Laplacian $-(-\Delta)^{\frac{s}{2}}u(x)$ can be approximated by
\begin{align}\label{FKGS:eq:2.8}
-(-\Delta)^{\frac{s}{2}}u_{N}(x_j)=-\sum\limits_{k=-N/2}^{N/2}|k\mu|^{s}\hat{u}(x_{j})e^{-\text{i}k\mu(x_{j}-x_{a})}.
\end{align}
 We denote ${u}_j={u}(x_j)$ and plug $\hat{u}_{k}$ in to \eqref{FKGS:eq:2.8}, then above approximation can be written as a matrix form. To this end,
\begin{align}\label{FKGS:eq:2.10}
-(-\Delta)^{\frac{s}{2}} u_{N}(x_j)&=-\sum\limits_{k=-N/2}^{N/2}|k\mu|^{s}(\frac{1}{Nc_{k}}\sum\limits_{l=0}^{N-1}{u_{l}}e^{-\text{i}k\mu(x-x_{a})})e^{-ik\mu(x_{j}-x_{a})}\nonumber\\
&=\sum\limits_{l=0}^{N-1}{u}_{l}(-\sum\limits_{k=-N/2}^{N/2}\frac{1}{Nc_{k}}|k\mu|^{s}e^{-\text{i}k\mu(x_j-x_{a})})=(D^{s}{U})_{j},
\end{align}
where the spectral differential matrix $D^{s}$ is an $N\times N$ symmetric matrix with the elements \cite{p1}
\begin{align*}
(D^{s})_{j,l}&=-\sum\limits_{k=-N/2}^{N/2}\frac{1}{Nc_{k}}|k\mu|^{\alpha}e^{-\text{i}k\mu(x_j-x_{l})}.
\end{align*}

\begin{rmk} \rm
When $s=2$, the matrix $D^{s}$ will reduce to the second-order spectral differential matrix \cite{p1}.
\end{rmk}

\subsection{Conservative semi-discrete scheme}

 We now discretize system \eqref{FKGS:eq:2.4}-\eqref{FKGS:eq:2.7} in space by using the Fourier pseudo-spectral method to construct a semi-discrete scheme as follows
\begin{align}\label{SG:eq:3.1}
U_t=2V,
\end{align}
\begin{align}\label{SG:eq:3.2}
V_t=\frac{1}{2}D^{\beta}U-\frac{1}{2}U+\frac{1}{2}(P^{2}+Q^{2}),
\end{align}
\begin{align}\label{SG:eq:3.3}
P_t=\frac{1}{2}D^{\alpha}Q+U\cdot Q,
\end{align}
\begin{align}\label{SG:eq:3.4}
Q_t=-\frac{1}{2}D^{\alpha}P-U\cdot P,
\end{align}
where $P^2=P\cdot P$, and `$\cdot$' means the point multiplication between vectors, that is, $U\cdot V=(u_{0}v_{0}, \cdots, u_{N-1}v_{N-1} )^{T}$. Let $Y=(U^{T},V^{T},P^{T},Q^{T})$, similar to the continuous case discussion, the above system can be rewritten as a canonical Hamiltonian form
\begin{align}\label{SG:eq:3.5}
\frac{dY}{dt}=f(Y)=S\nabla H(Y),
\end{align}
where  the discrete Hamiltonian is defined by
\begin{align}\label{SG:eq:3.6}
H(Y)=\frac{1}{4}\Big(-P^{T}D^{\alpha}P-U^{T}D^{\beta}U-Q^{T}D^{\alpha}Q+U^{T}U+4V^{T}V-2U^{T}(P^{2}+Q^{2})\Big),
\end{align}
and $S$ is a skew symmetric matrix with the following form
\begin{align*}
S=\left(\begin{array}{ccccc}
              {0}& I&0&0\\
              {-I}& 0&0&0 \\
              {0}& 0&0&-I\\
              {0}& 0&I&0
              \end{array}
\right).
\end{align*}
We define the mass
\begin{align}\label{SG:eq:3.7}
M(Y):=\|P\|^2+\|Q\|^2,
\end{align}
then the conservation laws of above semi-discrete scheme are given in the following theorem.

\begin{thm}\rm\label{2SG-lem2.1} The scheme \eqref{SG:eq:3.1}-\eqref{SG:eq:3.4} satisfies the semi-discrete mass and energy conservation law
\begin{align*}
\frac{dM(Y)}{dt}=0,\ \ \  \text{and} \ \  \frac{dH(Y)}{dt}=0.
\end{align*}

\end{thm}
\begin{prf}\rm
According to \eqref{SG:eq:3.7}, we can derive
\begin{align}\label{SG:eq:3.8}
\frac{dM(Y)}{dt}=2h\Big(\frac{1}{2}P^{T}D^{\alpha}Q-\frac{1}{2}Q^{T}D^{\alpha}P+P^{T}(U\cdot Q)-Q^{T}(P\cdot U)\Big)=0.
\end{align}
Then, based on the  skew symmetric of the matrix $S$, from \eqref{SG:eq:3.5}, we can obtain
\begin{align}\label{SG:eq:3.9}
\frac{dH(Y)}{dt}=\nabla H(Y)^{T}f(Y)=\nabla H(Y)^{T}S\nabla H(Y)=0,
\end{align}
We finish the proof.
\qed
\end{prf}

\subsection{Fully-discrete conservative schemes}

Based on the discussions of the semi-discrete scheme \eqref{SG:eq:3.1}-\eqref{SG:eq:3.4}, in this subsection, our goal is to establish a class of fully-discrete conservative schemes with regard to time and space.

In order to compare the standard AVF method and PAVF methods, we first use the original second-order AVF method to discrete the semi-discrete system \eqref{SG:eq:3.1}-\eqref{SG:eq:3.4} in time, and can obtain a fully-discrete scheme for the fractional KGS equation as follows
\begin{align}\label{SG:eq:3.10}
\delta_t U^m=2V^{m+\frac{1}{2}},
\end{align}
\begin{align}\label{SG:eq:3.11}
\delta_tV^m=\frac{1}{2}(D^{\beta}U^{m+\frac{1}{2}}-U^{m+\frac{1}{2}})&+\frac{1}{12}\Big((P^{m+1})^{2}+4(P^{m+\frac{1}{2}})^{2}+(P^{m})^2\nonumber\\
&+(Q^{m+1})^{2}+4(Q^{m+\frac{1}{2}})^{2}+(Q^{m})^2\Big),
\end{align}
\begin{align}\label{SG:eq:3.12}
\delta_t^{}P^m=\frac{1}{2}D^{\alpha}Q^{m+\frac{1}{2}}+\frac{1}{6}\Big(U^{m+1}\cdot Q^{n+1}+4U^{m+\frac{1}{2}}\cdot Q^{m+\frac{1}{2}}+U^{m}\cdot Q^{m}\Big),
\end{align}
\begin{align}\label{SG:eq:3.13}
\delta_t Q^m=-\frac{1}{2}D^{\alpha}P^{m+\frac{1}{2}}-\frac{1}{6}\Big(U^{m+1}\cdot P^{m+1}+4U^{m+\frac{1}{2}}\cdot P^{m+\frac{1}{2}}+U^{m}\cdot P^{m}\Big).
\end{align}
The proposed scheme \eqref{SG:eq:3.10}-\eqref{SG:eq:3.13} obtained by using the Fourier pseudo-spectral method to discrete the equation in space direction, and applying the
AVF method on the temporal direction, which is known as FAVF (Fourier pseudo-spectral AVF) scheme. It is observed that the scheme is a fully-implicit scheme and only can preserve the energy.

Next, we will present a fully-discrete scheme by applying the PAVF method to the semi-discrete system \eqref{SG:eq:3.1}-\eqref{SG:eq:3.4}. According to \eqref{FKGS:eq:3.4}, we have
\begin{align}\label{SG:eq:3.14}
\delta_t U^m=\int_{0}^{1}H_V(U^{m+1},\varepsilon V^{m+1}+(1-\varepsilon)V^{m}, P^{m},Q^{m})d\varepsilon,
\end{align}
\begin{align}\label{SG:eq:3.15}
\delta_t V^m=\int_{0}^{1}H_U(\varepsilon U^{m+1}+(1-\varepsilon)U^{m},V^{m}, P^{m},Q^{m})d\varepsilon,
\end{align}
\begin{align}\label{SG:eq:3.16}
\delta_t P^m=-\int_{0}^{1}H_Q( U^{m+1},V^{m+1}, P^{m+1},\varepsilon Q^{m+1}+(1-\varepsilon)Q^{m})d\varepsilon,
\end{align}
\begin{align}\label{SG:eq:3.17}
\delta_t Q^m=\int_{0}^{1}H_P( U^{m+1},V^{m+1}, \varepsilon P^{m+1}+(1-\varepsilon)P^{m},Q^{m})d\varepsilon,
\end{align}
which can be further integrated as
\begin{align}\label{SG:eq:3.18}
\delta_tU^m=2V^{m+\frac{1}{2}},
\end{align}
\begin{align}\label{SG:eq:3.19}
\delta_tV^m=\frac{1}{2}\Big(D^{\beta}U^{m+\frac{1}{2}}-U^{m+\frac{1}{2}}+(P^{m})^2+(Q^{m})^2),
\end{align}
\begin{align}\label{SG:eq:3.20}
\delta_tP^m=\frac{1}{2}D^{\alpha}Q^{m+\frac{1}{2}}+U^{m+1}\cdot Q^{m+\frac{1}{2}},
\end{align}
\begin{align}\label{SG:eq:3.21}
\delta_tQ^m=-\frac{1}{2}D^{\alpha}P^{m+\frac{1}{2}}-U^{m+1}\cdot P^{m+\frac{1}{2}}.
\end{align}
The proposed scheme \eqref{SG:eq:3.18}-\eqref{SG:eq:3.21} is
 known as FPAVF (Fourier pseudo-spectral PAVF) scheme. A result on mass and energy of conservation properties  of the scheme is presented by the following theorem.

\begin{thm}\rm \label{FSE:thm3.2} The FPAVF scheme \eqref{SG:eq:3.18}-\eqref{SG:eq:3.21} possesses the following discrete total mass and energy conservation law
\begin{align*}
{{M}}^m={{M}}^{m+1},\ \ \ H^m=H^{m+1}, \ \ 0\le m\le M-1,
\end{align*}
where the mass is defined as
 \begin{align*}
 {{M}}^m=\|P^m\|^2+\|Q^m\|^2,
\end{align*}
and the energy has the form
 \begin{align*}
 {H}^m=\frac{1}{4}\Big(&-(P^{m})^{T}D^{\alpha}P^{m}-(U^{m})^{T}D^{\beta}U^{m}-(Q^{m})^{T}D^{\alpha}Q^{m}\nonumber\\
&+(U^{m})^{T}U^{m}+4(V^{m})^{T}V^{m}-2(U^{m})^{T}\big((P^{m})^{2}+(Q^{m})^{2}\big)\Big).
\end{align*}
\end{thm}
\begin{prf}\rm Taking the inner product of \eqref{SG:eq:3.20}-\eqref{SG:eq:3.21} with $P^{m+1}+P^m$ and $Q^{m+1}+Q^m$, respectively, we have
\begin{align}\label{NLS:eq:3.22}
\frac{\tau}{h}(P^{m+1}+P^m)(P^{m+1}-P^m)=h(P^{m+1}+P^m)^T(\frac{1}{2}D^{\alpha}Q^{m+\frac{1}{2}}+U^{m+1}\cdot Q^{m+\frac{1}{2}}),
\end{align}
and
\begin{align}\label{NLS:eq:3.23}
\frac{\tau}{h}(Q^{m+1}+Q^m)(Q^{m+1}-Q^m)=h(Q^{m+1}+Q^m)^T(-\frac{1}{2}D^{\alpha}P^{m+\frac{1}{2}}-U^{m+1}\cdot P^{m+\frac{1}{2}}).
\end{align}
Adding them together, we deduce
\begin{align}\label{NLS:eq:3.24}
\frac{1}{\tau}(\|P^{m+1}\|^2+\|Q^{m+1}\|^2-\|P^{m}\|^2-\|Q^{m}\|^2)=0,
\end{align}
which implies that
\begin{align*}
 {{M}}^m={{M}}^{m+1}.
\end{align*}
Then, we take the discrete inner products of \eqref{SG:eq:3.18}-\eqref{SG:eq:3.21} with $\delta_t V^m$, $\delta_t U^m$, $\delta_t Q^m$ and $-\delta_t P^m$, respectively, and summing them together, we can obtain the energy conservation law
\begin{align*}
 {H}^m={H}^{m+1}.
\end{align*}
This ends the proof.
\qed
\end{prf}

 To derive the second-order temporal accuracy conservative PAVF schemes for the fractional KGS equation, we first propose the adjoint of the FPAVF scheme \eqref{SG:eq:3.18}-\eqref{SG:eq:3.21} as follows
\begin{align}\label{SG:eq:3.25}
\delta_t U^m=2V^{m+\frac{1}{2}},
\end{align}
\begin{align}\label{SG:eq:3.26}
\delta_t V^m=\frac{1}{2}\Big(D^{\beta}U^{m+\frac{1}{2}}-U^{m+\frac{1}{2}}+(P^{m+1})^2+(Q^{m+1})^2\Big),
\end{align}
\begin{align}\label{SG:eq:3.27}
\delta_t P^m=\frac{1}{2}D^{\alpha}Q^{m+\frac{1}{2}}+U^{m}\cdot Q^{m+\frac{1}{2}},
\end{align}
\begin{align}\label{SG:eq:3.28}
\delta_t Q^m=-\frac{1}{2}D^{\alpha}P^{m+\frac{1}{2}}-U^{m}\cdot P^{m+\frac{1}{2}}.
\end{align}
Then, combing the FPAVF scheme and the adjoint FPAVF scheme, the FPAVF-C (Fourier pseudo-spectral PAVF-C) scheme can be given by
\begin{align}\label{SG:eq:3.56}
&\frac{1}{\tau}(U^{*}-U^{m})=V^{*}+V^{m},\\
&\frac{1}{\tau}(V^{*}-V^{m})=\frac{1}{4}\Big(D^{\beta}(U^{*}+U^{m})-(U^{*}+U^{m})+2(P^{m})^2+2(Q^{m})^2\Big),\\
&\frac{1}{\tau}(P^{*}-P^{m})=\frac{1}{4}D^{\alpha}(Q^{*}+Q^{m})+\frac{1}{2}U^{*}\cdot (Q^{*}+Q^{m}),\\
&\frac{1}{\tau}(Q^{*}-Q^{m})=-\frac{1}{4}D^{\alpha}(P^{*}+P^{m})-\frac{1}{2}U^{*}\cdot (P^{*}+P^{m}),\\
&\frac{1}{\tau}(U^{m+1}-U^{*})=V^{*}+V^{m+1},\\
&\frac{1}{\tau}(V^{m+1}-V^{*})=\frac{1}{4}\Big(D^{\beta}(U^{*}+U^{m+1})-(U^{*}+U^{m+1})+2(P^{m+1})^2+2(Q^{m+1})^2\Big),\\
&\frac{1}{\tau}(P^{m+1}-P^{*})=\frac{1}{4}D^{\alpha}(Q^{*}+Q^{m+1})+\frac{1}{2}U^{*}\cdot (Q^{*}+Q^{m+1}),\\
&\frac{1}{\tau}(Q^{m+1}-Q^{*})=-\frac{1}{4}D^{\alpha}(P^{*}+P^{m+1})-\frac{1}{2}U^{*}\cdot (P^{*}+P^{m+1}),
\end{align}
and the FAVF-P (Fourier pseudo-spectral PAVF-P) scheme
\begin{align}\label{SG:eq:3.47}
\delta_tU^m=V^{m+\frac{1}{2}},
\end{align}
\begin{align}\label{SG:eq:3.48}
\delta_t V^m=\frac{1}{2}(D^{\beta}U^{m+\frac{1}{2}}-U^{m+\frac{1}{2}})+\frac{1}{4}\Big((P^{m})^2+(Q^{m})^2+(P^{m+1})^2+(Q^{m+1})^2),
\end{align}\begin{align}\label{SG:eq:3.49}
\delta_t P^m=\frac{1}{2}D^{\alpha}Q^{m+\frac{1}{2}}+U^{m+\frac{1}{2}}\cdot Q^{m+\frac{1}{2}},
\end{align}
\begin{align}\label{SG:eq:3.50}
\delta_t Q^m=-\frac{1}{2}D^{\alpha}P^{m+\frac{1}{2}}-U^{m+\frac{1}{2}}\cdot P^{m+\frac{1}{2}}.
\end{align}
It is observe that the FPAVF-C scheme and the FPAVF-P scheme both are symmetric methods of second-order accuracy and can preserve the discrete mass and energy conservation law.

\begin{rmk} \rm
The FAVF scheme does not preserve the mass conservation law, the FPAVF scheme, including the FPAVF-C scheme and the FPAVF-P scheme, can not only preserve the energy but also preserve the mass.

\end{rmk}

\section{Numerical examples}
In this section,  numerical examples of the conservative schemes are given to support our theoretical analysis.
 Fast solver method based on the fast Fourier transformation (FFT) technique is used in practical computation, which can reduce the memory requirement and the computational complexity.
To obtain numerical errors, we use the error function defined as follows
\begin{align*}
&E(\tau)=\|U_{N}^{M}-U_{N}^{2M}\|_{\infty}+\|P_{N}^{M}-P_{N}^{2M}\|_{\infty}+\|Q_{N}^{M}-Q_{N}^{2M}\|_{\infty},\\
&E(N)=\|U^{M}_{N}-U^{M}_{2N}\|_{\infty}+\|P^{M}_{N}-P^{M}_{2N}\|_{\infty}+\|Q^{M}_{N}-Q^{M}_{2N}\|_{\infty},
\end{align*}
where $\|U_{N}^{M}-U_{N}^{2M}\|_{\infty}:=\|U(\frac{T}{M},\frac{L}{N})-U(\frac{T}{2M},\frac{L}{N})\|_{\infty}$, $\|U^{M}_{N}-U^{M}_{2N}\|_{\infty}:=\|U(\frac{T}{M},\frac{L}{N})-U(\frac{T}{M},\frac{L}{2N})\|_{\infty}$, etc.,
and the convergence orders in time and space of the $l^{\infty}$-norm errors on two successive time step sizes $\tau$ and $\tau/2$ and two successive polynomial degrees $N$ and $2N$ are calculated as
\begin{align*}
\text{order}=\left\{\begin{array}{lll}
              &{\text{log}_{2}{[{E}(\tau)/{E}(\tau/2)]}},\ &\text{ in time}, \\
 \\
&{\text{log}_{2}[{E}( N)/{E}( 2N)]},\ &\text{ in space}.
\\
             \end{array}
\right.
\end{align*}
The relative errors of energy and mass
are defined as
\begin{align*}
RH^{m}=|(H^{m}-H^{0})/H^{0}|,\ \ \ RM^{m}=|(M^{m}-M^{0})/M^{0}|,
\end{align*}
where $H^{m}$ and $M^{m}$ denote the energy and mass at $t=m\tau$, respectively.

\textbf{Example 4.1.}
 We study one-dimensional case of system \eqref{FKGS:eq:1.1}-\eqref{FKGS:eq:1.2} with the initial conditions are chosen as
\begin{align*}
&\varphi(x,0)=\frac{3\sqrt{2}}{4\sqrt{1-r^2}}\text{sech}^{2}\Big(\frac{1}{2\sqrt{1-r^2}}(x-x_{0})\Big)\text{exp}(\text{i}rx),\\
&u(x,0)=\frac{3r}{4(1-r^2)}\text{sech}^{2}\Big(\frac{1}{2\sqrt{1-r^2}}(x-x_{0})\Big),\\
&u_{t}(x,0)=\frac{3}{4(1-r^2)^{3/2}}\text{sech}^{2}\Big(\frac{1}{2\sqrt{1-r^2}}(x-x_{0})\text{tanh}(\frac{1}{2\sqrt{1-r^2}}(x-x_{0})\Big),
\end{align*}
where $|r|<1$ is the propagating velocity of the wave and $x_0$ is the initial phase. In our computation, we take the computation domain $\Omega=(-20, 20)$. Without loss of generality, we take $\alpha=\beta=1.4, 1.7, 2$.
If $\alpha=2$, system \eqref{FKGS:eq:1.1}-\eqref{FKGS:eq:1.2} reduces to the classical KGS equation with the exact solutions
\begin{align*}
&\varphi(x,t)=\frac{3\sqrt{2}}{4\sqrt{1-r^2}}\text{sech}^{2}\Big(\frac{1}{2\sqrt{1-r^2}}(x-rt-x_{0})\Big)\text{exp}\Big(\text{i}(rx+\frac{1-r^2+r^4}{2(1-r^2)}t)\Big),\\
&u(x,t)=\frac{3}{4(1-r^2)}\text{sech}^{2}\Big(\frac{1}{2\sqrt{1-r^2}}(x-rt-x_{0})\Big).
\end{align*}

First, to test the accuracy of the our scheme, we consider initial value $r=-0.8,~x_0=0$. Fig. 1 plots the time errors of the four conservative schemes in the log scale. It is demonstrated that the FPAVF scheme for different $\alpha$ is first-order of convergence in time, and others schemes are second-order accuracy in time. Then we test the spatial accuracy of the four methods. Here we set the $\tau=10^{-6}$ and plot the spatial errors at $T=1$ for different $\alpha$ in Fig. 2,  which implies that for the sufficiently smooth problem, our conservative schemes are of spectral accuracy in spatial.
To compare the efficiency of the four conservative schemes, in Fig. 3, we present comparisons on the computational costs of the schemes with different time step. As expected, the FPAVF scheme requires the smallest computational cost, regardless of its first-order accuracy. Among the second-order ones, the FPAVF-C scheme is much more efficient than the FPAVF-P scheme and the most time-consuming FPAVF scheme. As a consequence, the FPAVF-C scheme is obviously the optimal choice among the four energy-preserving methods which will exhibit greater advantages of lower computational cost in simulating multi-dimensional problems.

%Obviously, we can conclude that the FPAVF scheme, including the FPAVF-C scheme and FPAVF-P scheme can reduce the computational cost.
%This is due to the AVF method and the PAVF-P method are fully implicit, their computational costs are comparable but more than the PAVF method and the PAVF-C method which are only linearly implicit.

\begin{figure}[H]
\centering
\begin{minipage}[t]{60mm}
\includegraphics[width=60mm]{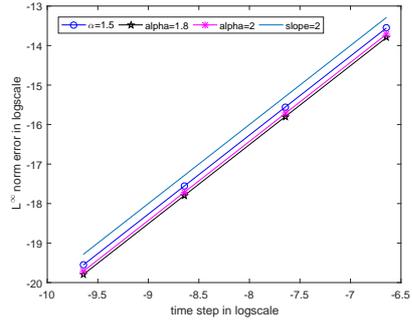}\\
{\footnotesize  \centerline {(a) FAVF scheme}}
\end{minipage}
\begin{minipage}[t]{60mm}
\includegraphics[width=60mm]{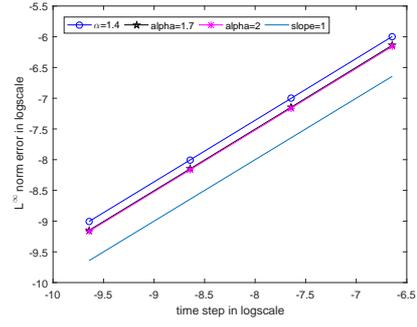}\\
{\footnotesize  \centerline {(b) FPAVF scheme}}
\end{minipage}
\begin{minipage}[t]{60mm}
\includegraphics[width=60mm]{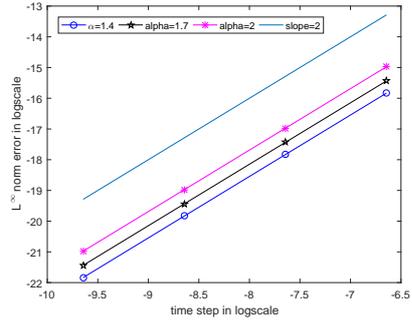}\\
{\footnotesize  \centerline {(c) FAVF-C scheme}}
\end{minipage}
\begin{minipage}[t]{60mm}
\includegraphics[width=60mm]{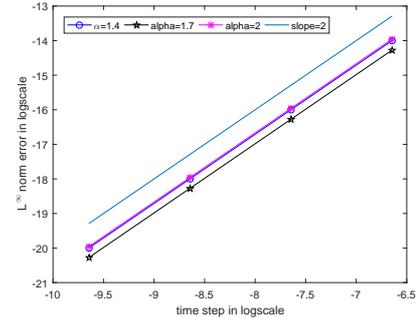}\\
{\footnotesize  \centerline {(d) FAVF-P scheme}}
\end{minipage}
\caption{\small {\small{Temporal accuracy of four schemes for different $\alpha$ with $N=128$.} }}\label{fig522}
\end{figure}

\begin{figure}[H]
\centering
\begin{minipage}[t]{60mm}
\includegraphics[width=60mm]{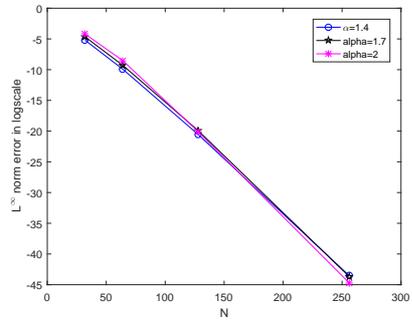}\\
{\footnotesize  \centerline {(a) FAVF scheme}}
\end{minipage}
\begin{minipage}[t]{60mm}
\includegraphics[width=60mm]{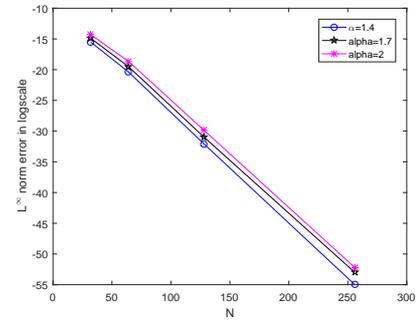}\\
{\footnotesize  \centerline {(b) FPAVF scheme}}
\end{minipage}
\begin{minipage}[t]{60mm}
\includegraphics[width=60mm]{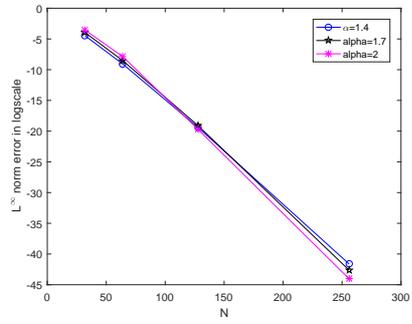}\\
{\footnotesize  \centerline {(c) FAVF-C scheme}}
\end{minipage}
\begin{minipage}[t]{60mm}
\includegraphics[width=60mm]{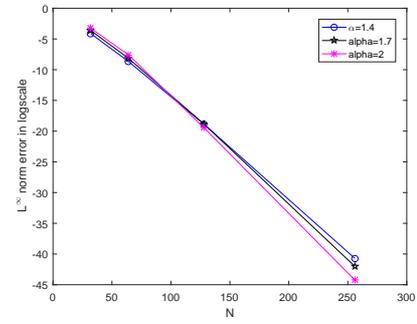}\\
{\footnotesize  \centerline {(d) FAVF-P scheme}}
\end{minipage}
\caption{\small {\small{Spatial accuracy of four schemes for different $\alpha$ with $\tau=10^{-6}$.} }}\label{fig522}
\end{figure}

\begin{figure}[H]
\centering\begin{minipage}[t]{60mm}
\includegraphics[width=60mm]{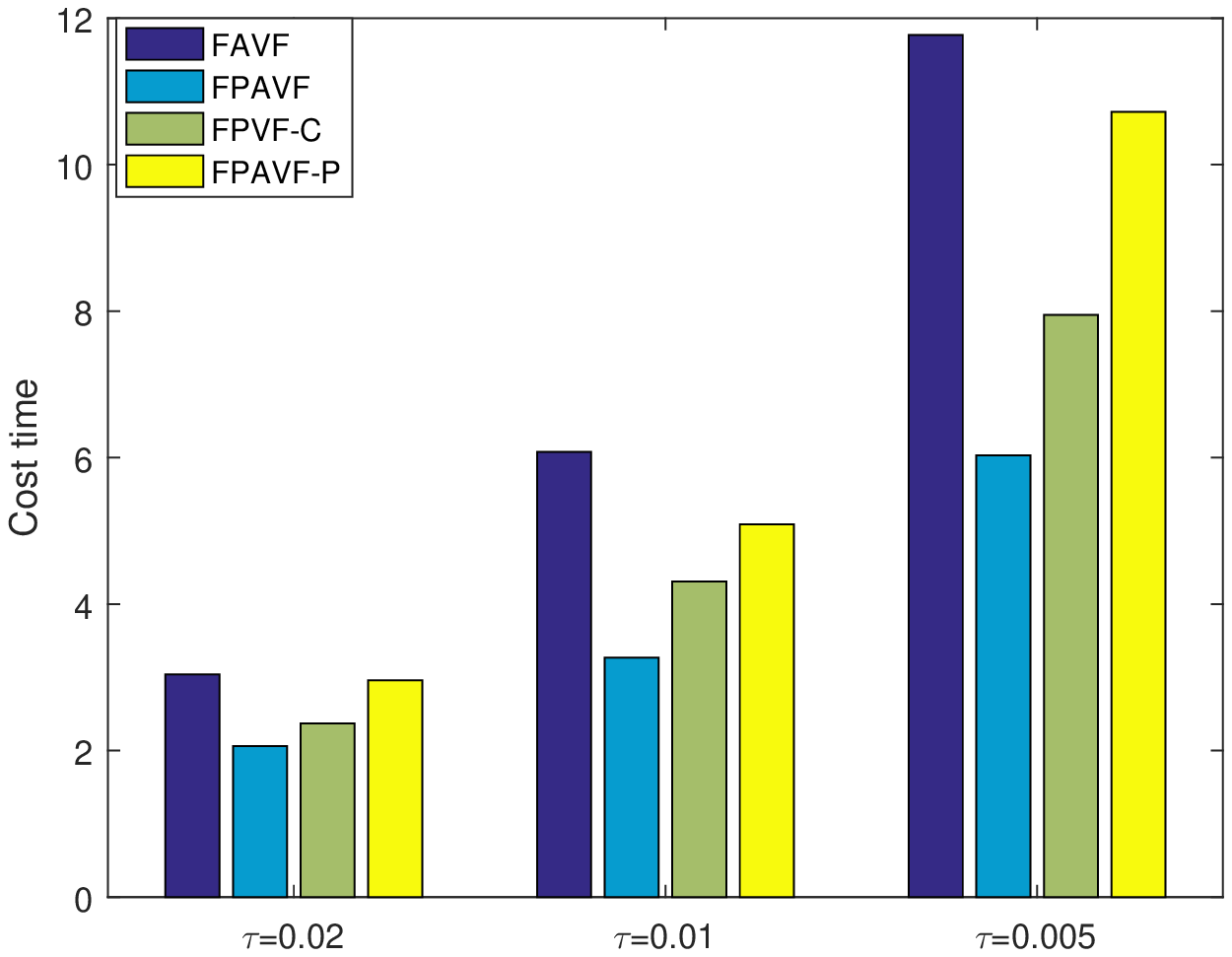}\\
{\footnotesize  \centerline {(a) $\alpha=1.7$}}
\end{minipage}
\begin{minipage}[t]{60mm}
\includegraphics[width=60mm]{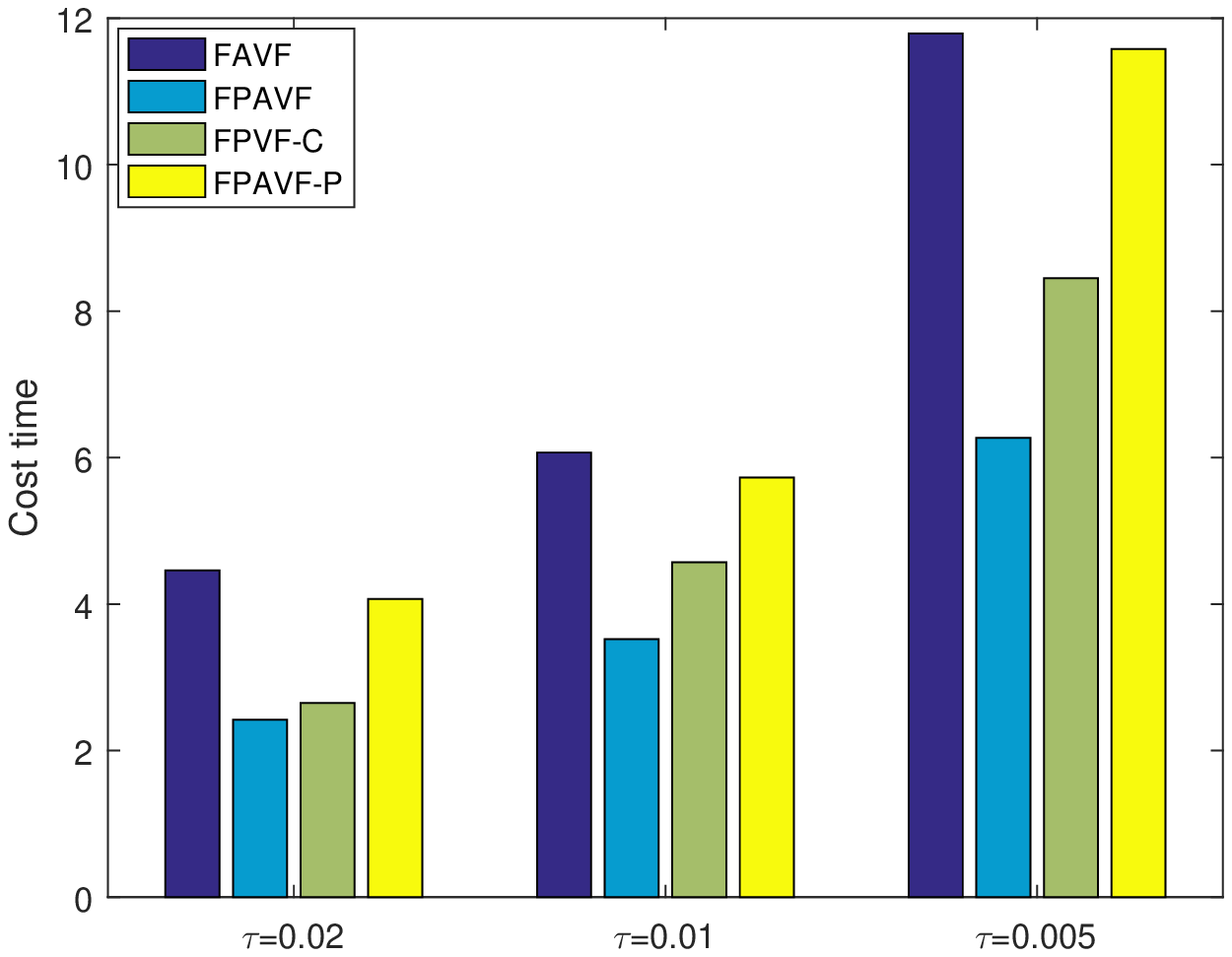}\\
{\footnotesize  \centerline {(b) $\alpha=2$}}
\end{minipage}
\caption{\small {\small{Computational cost of four schemes with different time step till $T=50$.} }}\label{fig522}
\end{figure}

Second, we choose $T=50$ and verify the discrete conservation laws of the proposed schemes. We take $N=128,$~$\tau=0.001$ and compute the discrete mass and energy. Fig. 4 and Fig. 5 show the relative errors of mass $M$ and energy $H$ for different $\alpha$, respectively.
The pictures demonstrate that all schemes can preserve the energy very well in discrete sense. However, the FAVF scheme fails to preserve the discrete mass conservation law, in contrast, the rest three schemes can preserve the mass exactly.

\begin{figure}[H]
\centering
\begin{minipage}[t]{50mm}
\includegraphics[width=50mm]{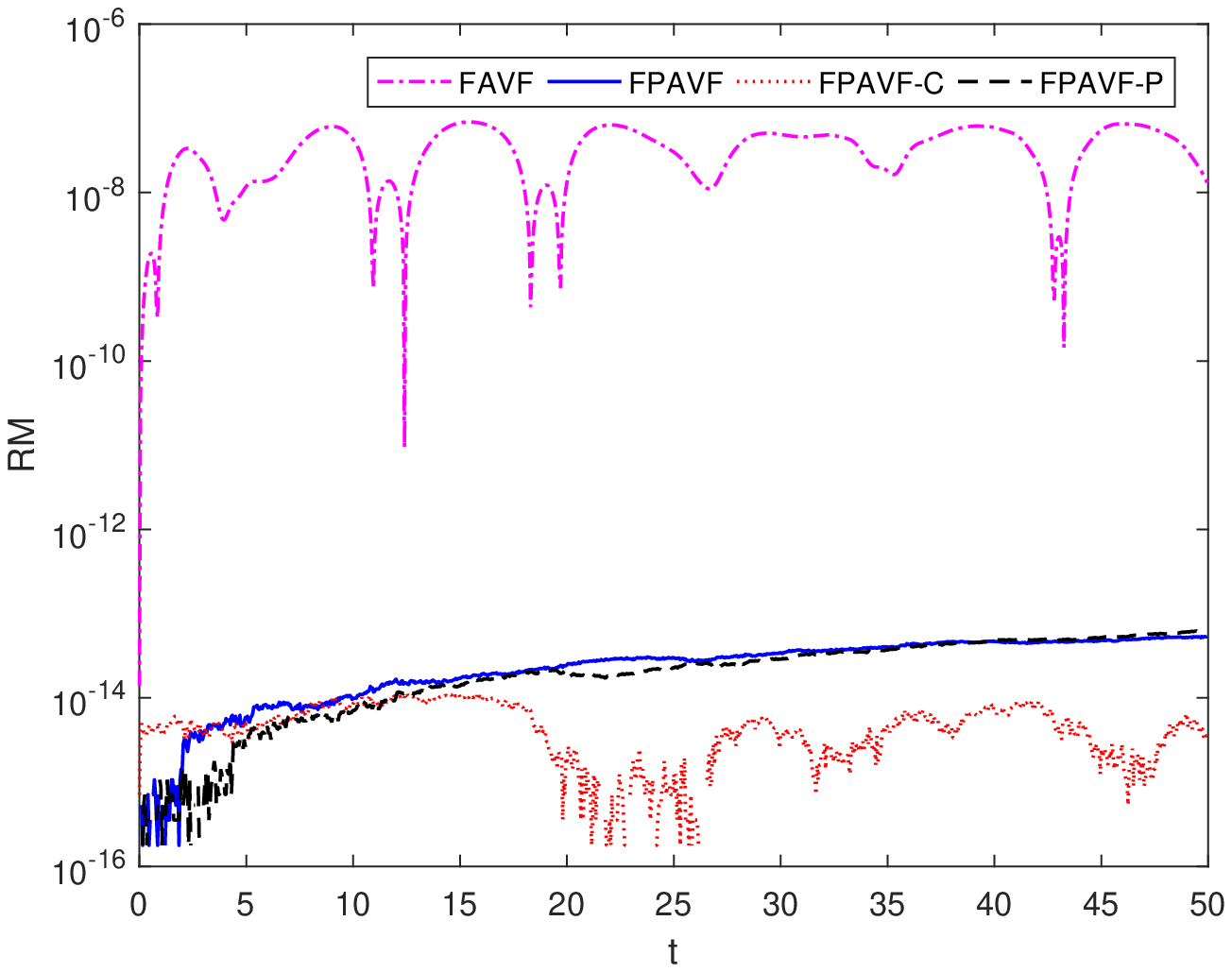}\\
{\footnotesize  \centerline {(a) $\alpha=1.4$}}
\end{minipage}
\begin{minipage}[t]{50mm}
\includegraphics[width=50mm]{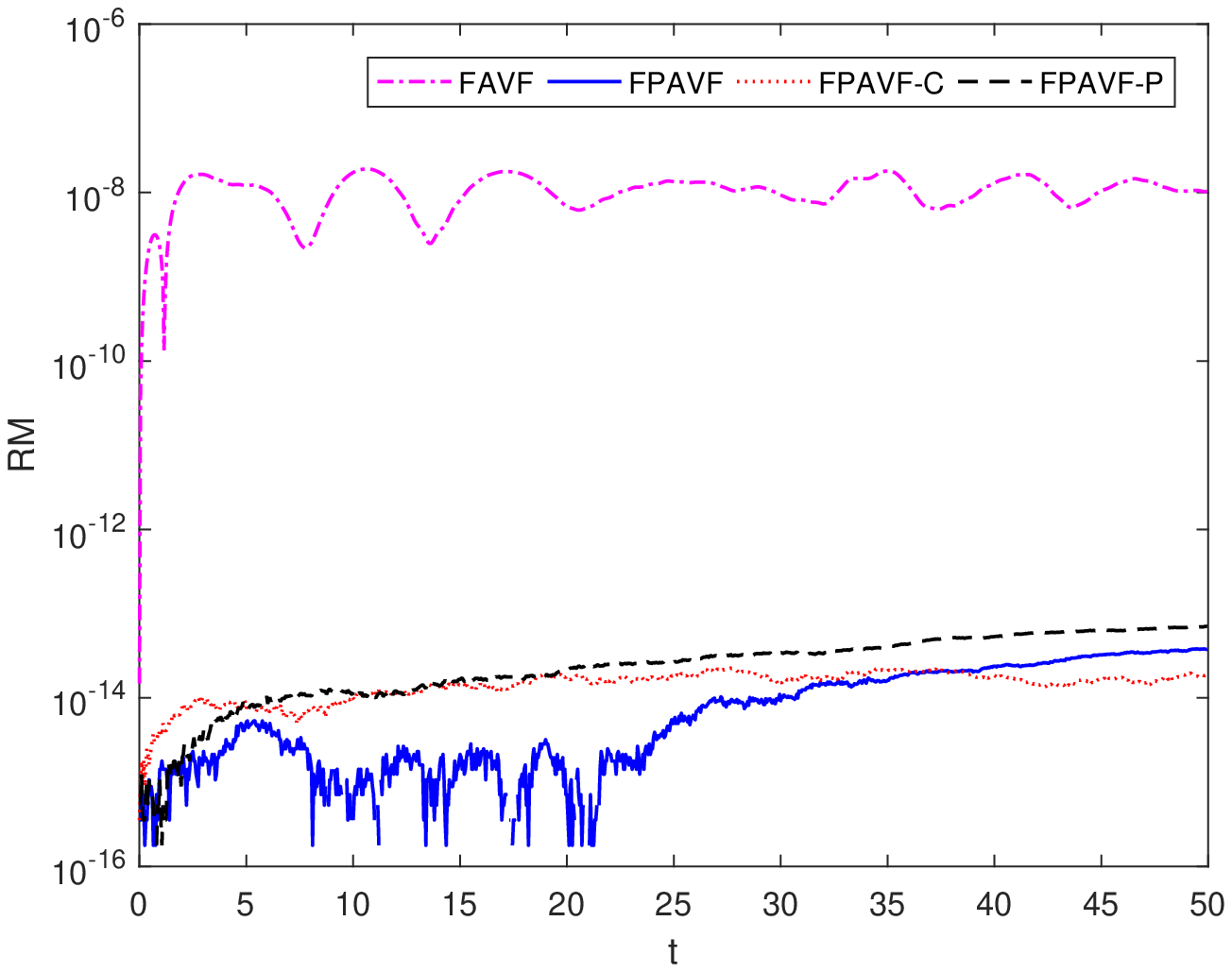}\\
{\footnotesize  \centerline {(b) $\alpha=1.7$}}
\end{minipage}
\begin{minipage}[t]{50mm}
\includegraphics[width=50mm]{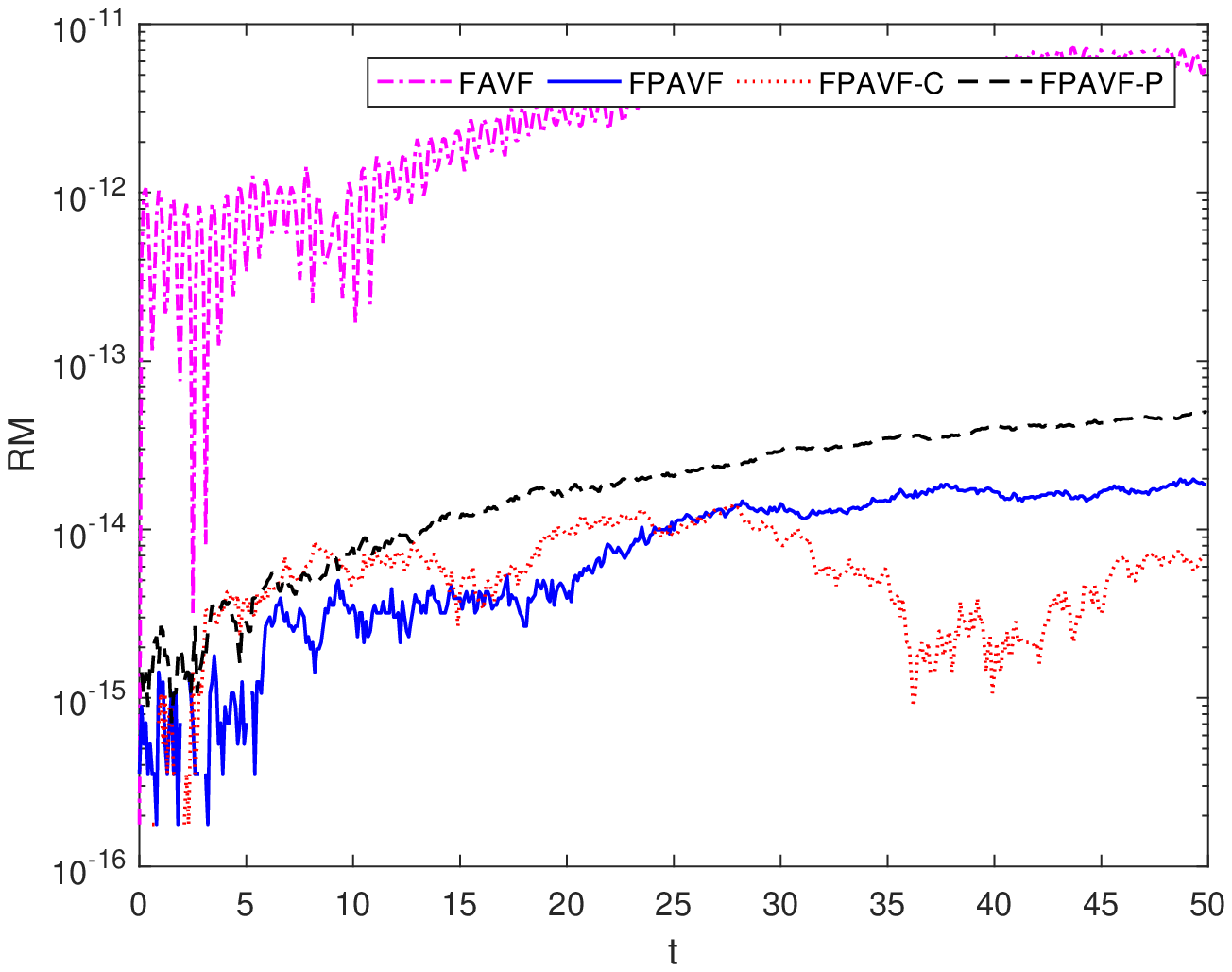}\\
{\footnotesize  \centerline {(c) $\alpha=2$}}
\end{minipage}
\caption{\small {The relative errors of mass with $N=128,~\tau=0.001$ for different $\alpha$. }}\label{fig522}
\end{figure}

\begin{figure}[H]
\centering
\begin{minipage}[t]{50mm}
\includegraphics[width=50mm]{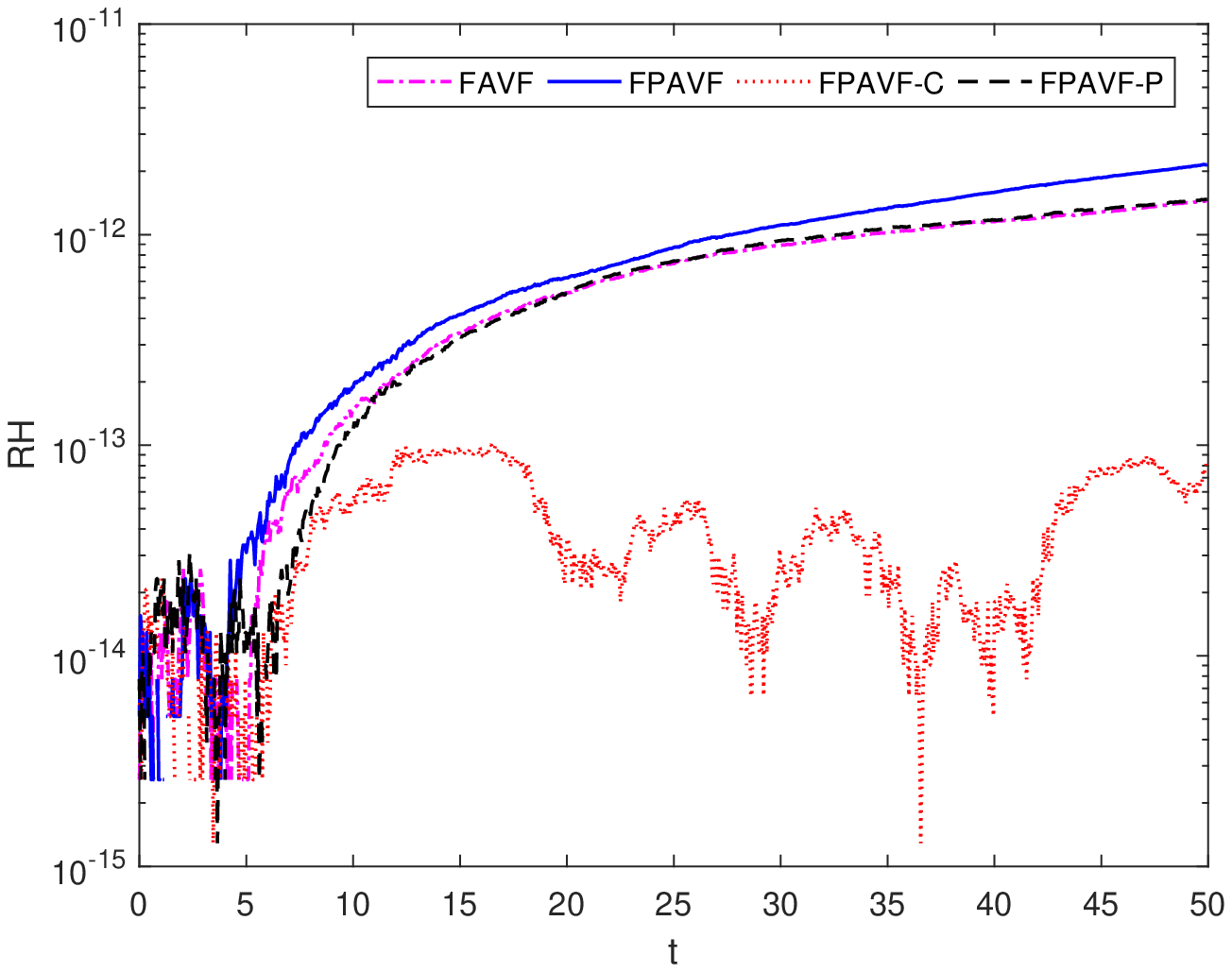}\\
{\footnotesize  \centerline {(a) $\alpha=1.4$}}
\end{minipage}
\begin{minipage}[t]{50mm}
\includegraphics[width=50mm]{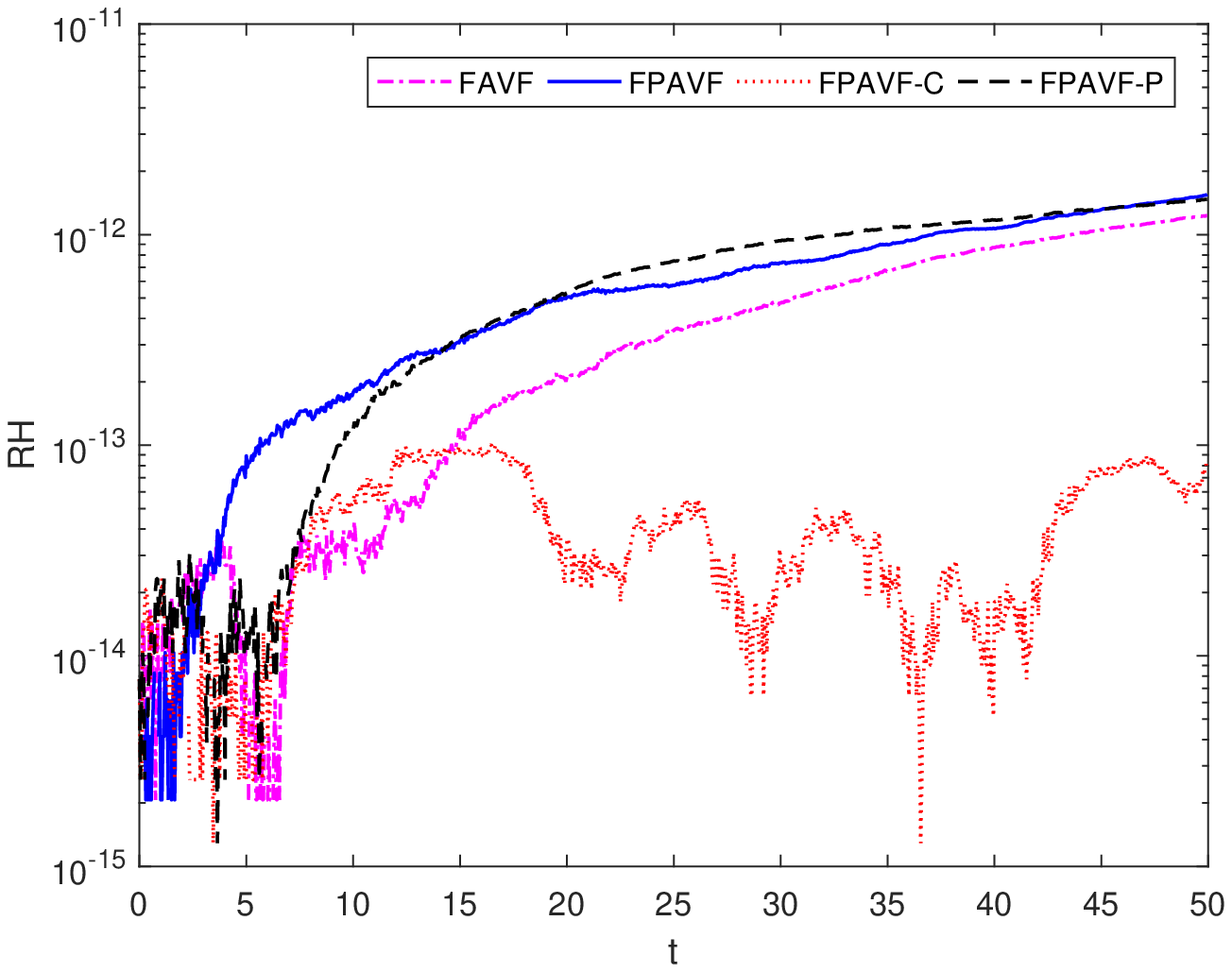}\\
{\footnotesize  \centerline {(b) $\alpha=1.7$}}
\end{minipage}
\begin{minipage}[t]{50mm}
\includegraphics[width=50mm]{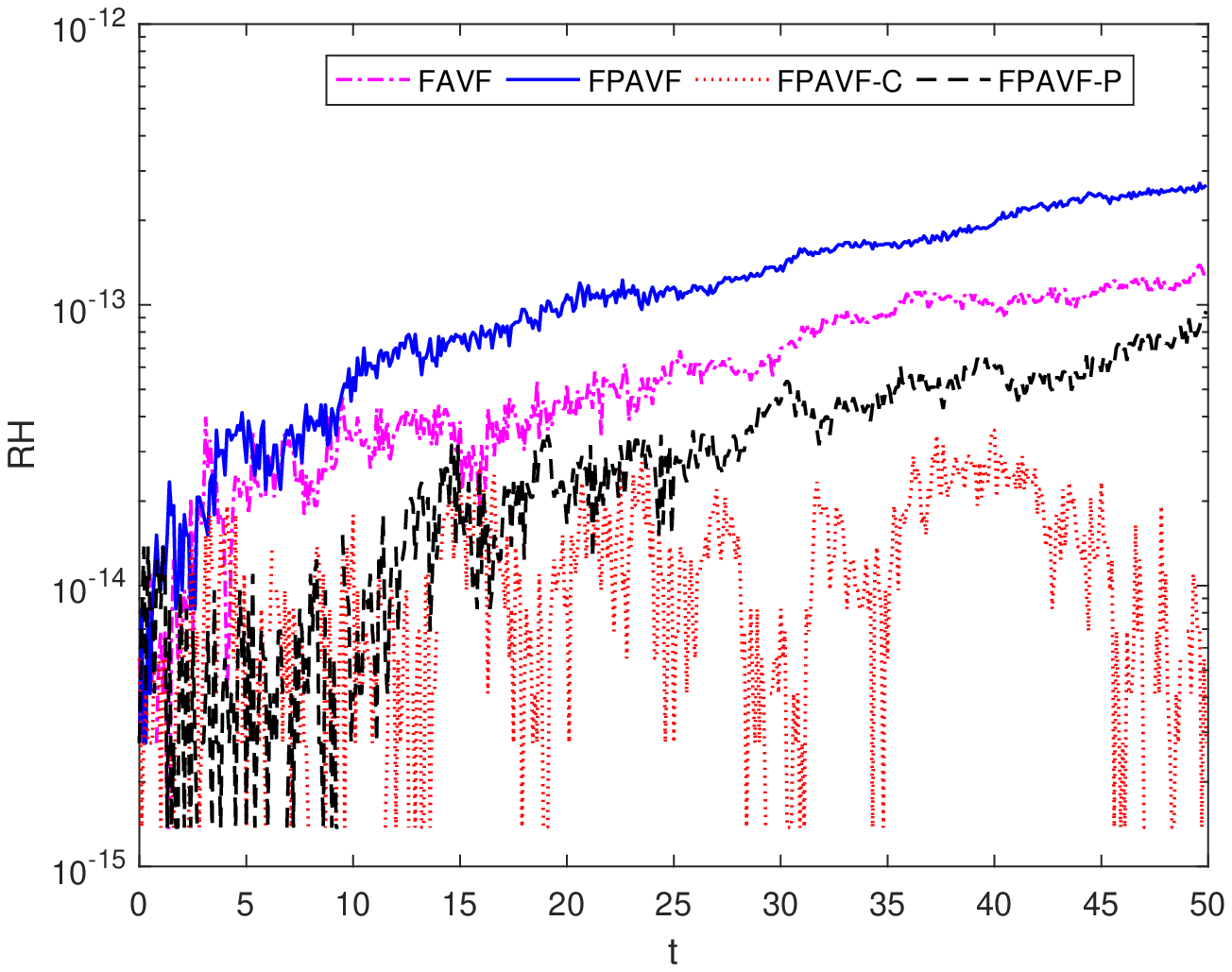}\\
{\footnotesize  \centerline {(c) $\alpha=2$}}
\end{minipage}
\caption{\small {The relative errors of energy with $N=128,~\tau=0.001$ for different $\alpha$. }}\label{fig522}
\end{figure}

Finally, we investigate the relationship between the fractional order $\alpha$ and the shape of the solition for the problem.  Here we take $T=10$. The wave profiles of $|\varphi^n|=|Q^n+\text{i}P^n|$ and $U^n$ obtained by four schemes are presented in Fig. 6, respectively, it demonstrates that the fractional order $\alpha$ will affect the shape of the soliton, and the shape of the soliton will change more quickly when $\alpha$ becomes smaller. These properties are similar as the numerical simulations of the fractional Schr\"{o}dinger equation \cite{p11,p12}, which can be used in physics to modify the shape of wave without
change of the nonlinearity and dispersion effects.

\begin{figure}[H]
\centering
\begin{minipage}[t]{60mm}
\includegraphics[width=60mm]{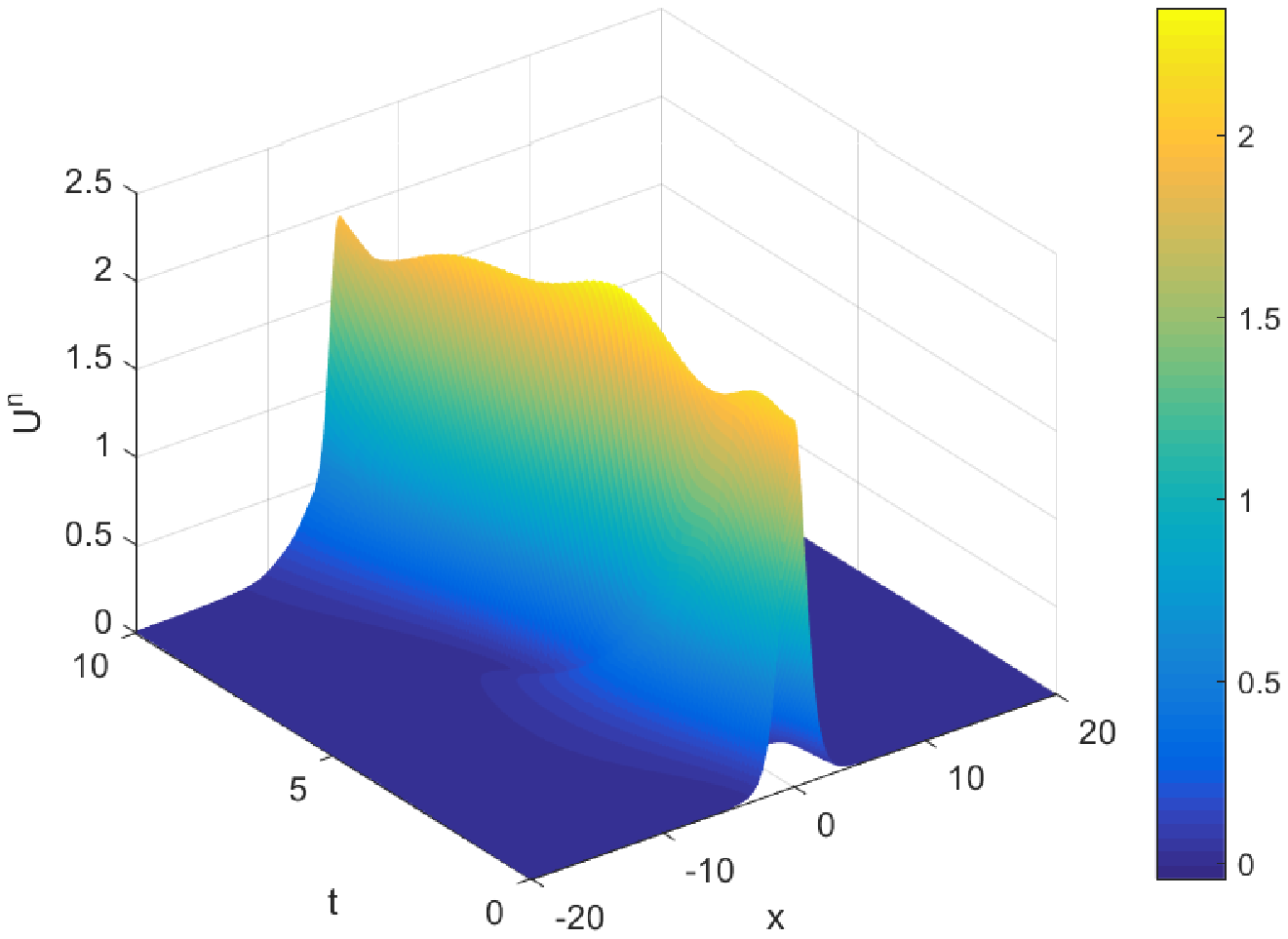}\\
{\footnotesize  \centerline {(a) $\alpha=1.4$}}
\end{minipage}
\begin{minipage}[t]{60mm}
\includegraphics[width=60mm]{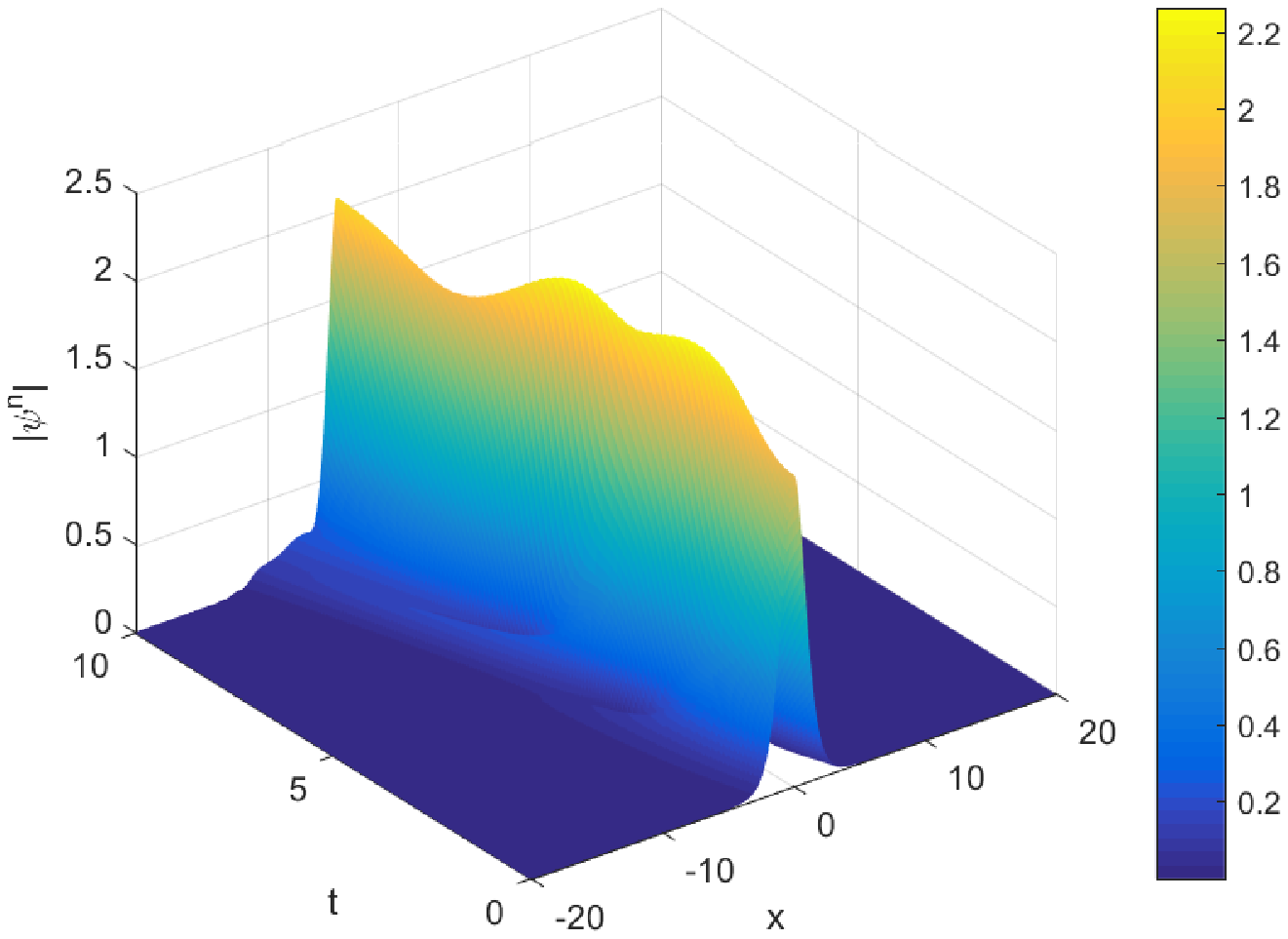}\\
{\footnotesize  \centerline {(b) $\alpha=1.4$}}
\end{minipage}
\begin{minipage}[t]{60mm}
\includegraphics[width=60mm]{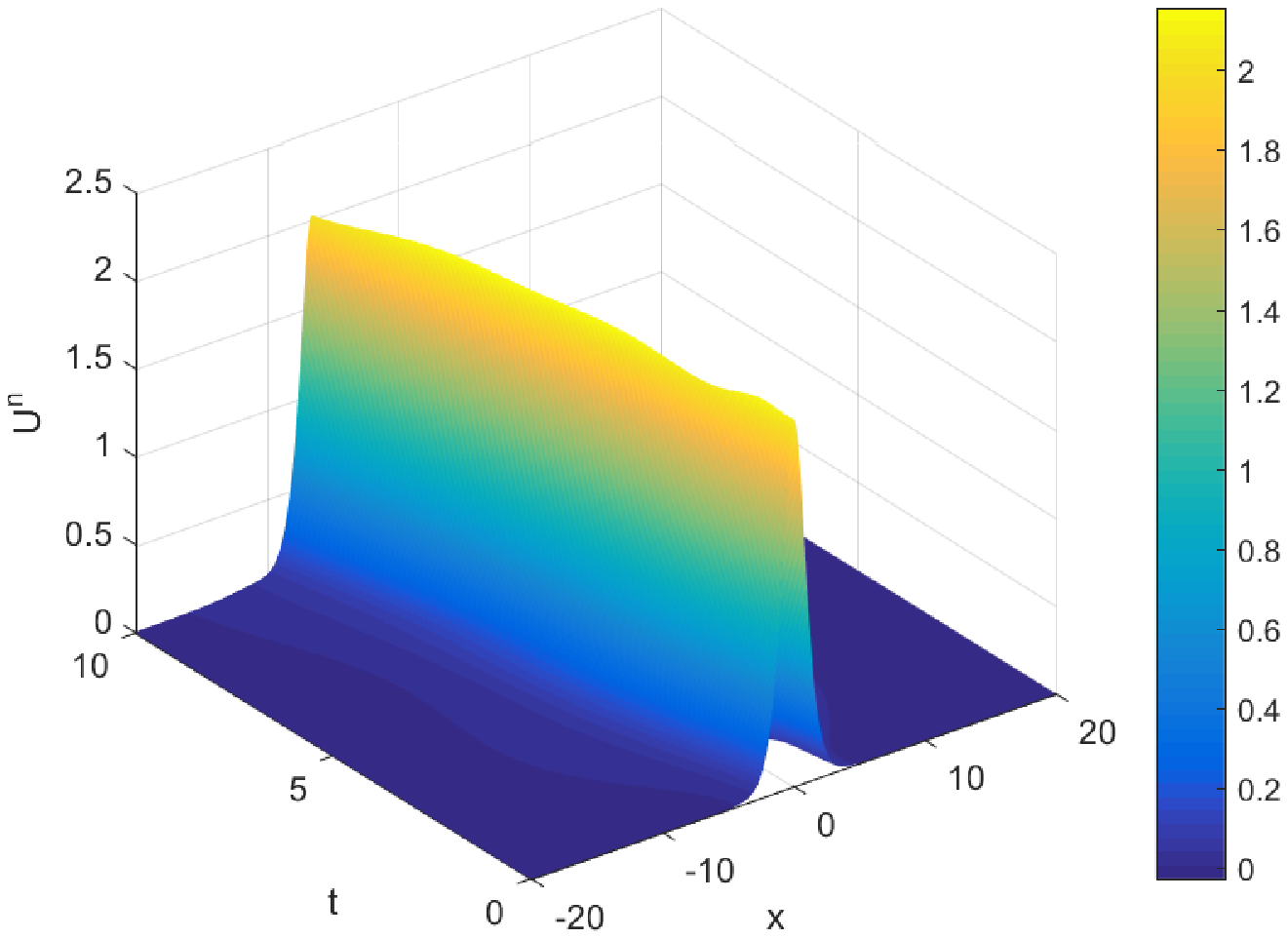}\\
{\footnotesize  \centerline {(c) $\alpha=1.7$}}
\end{minipage}
\begin{minipage}[t]{60mm}
\includegraphics[width=60mm]{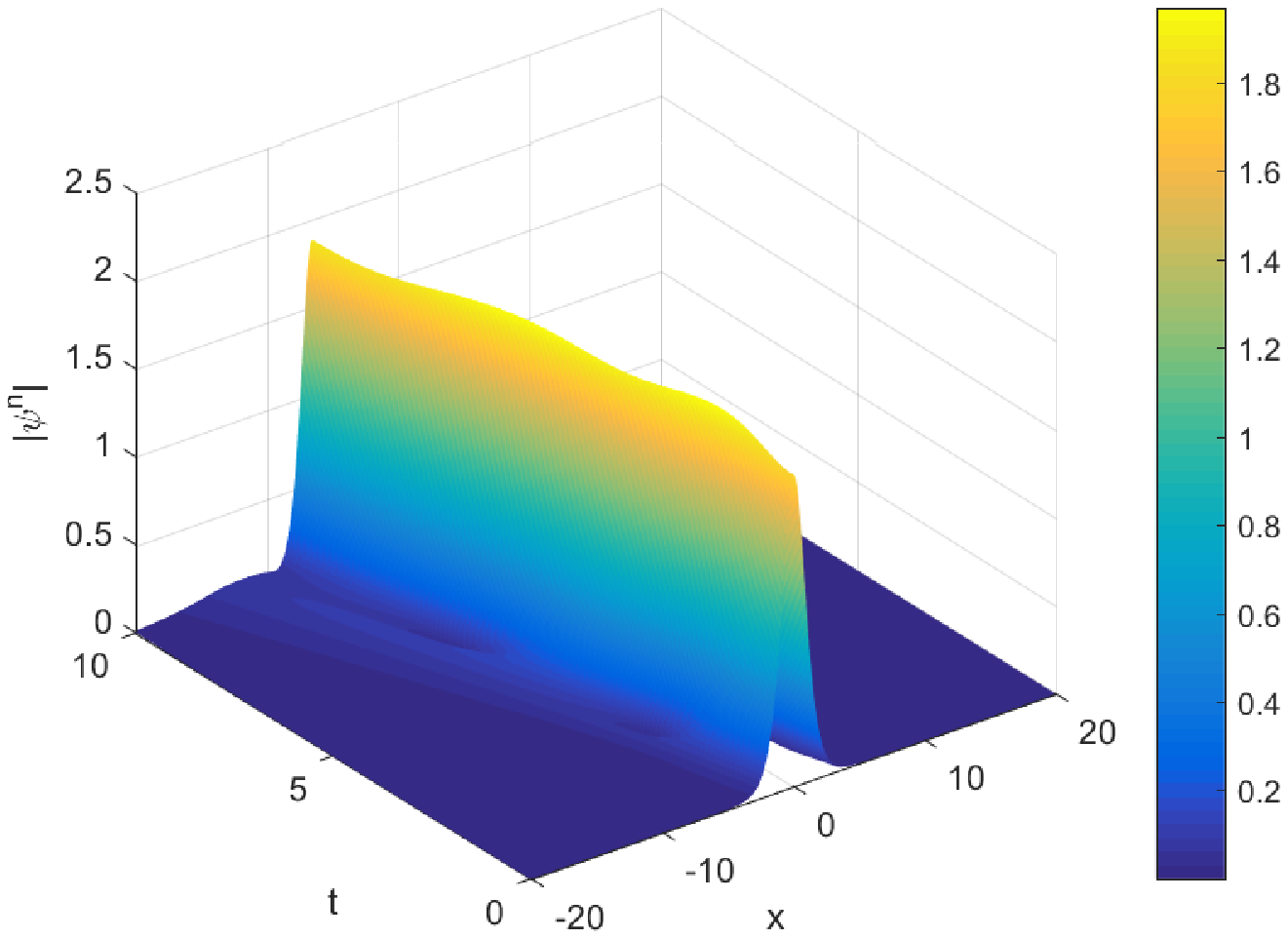}\\
{\footnotesize  \centerline {(d) $\alpha=1.7$}}
\end{minipage}
\begin{minipage}[t]{60mm}
\includegraphics[width=60mm]{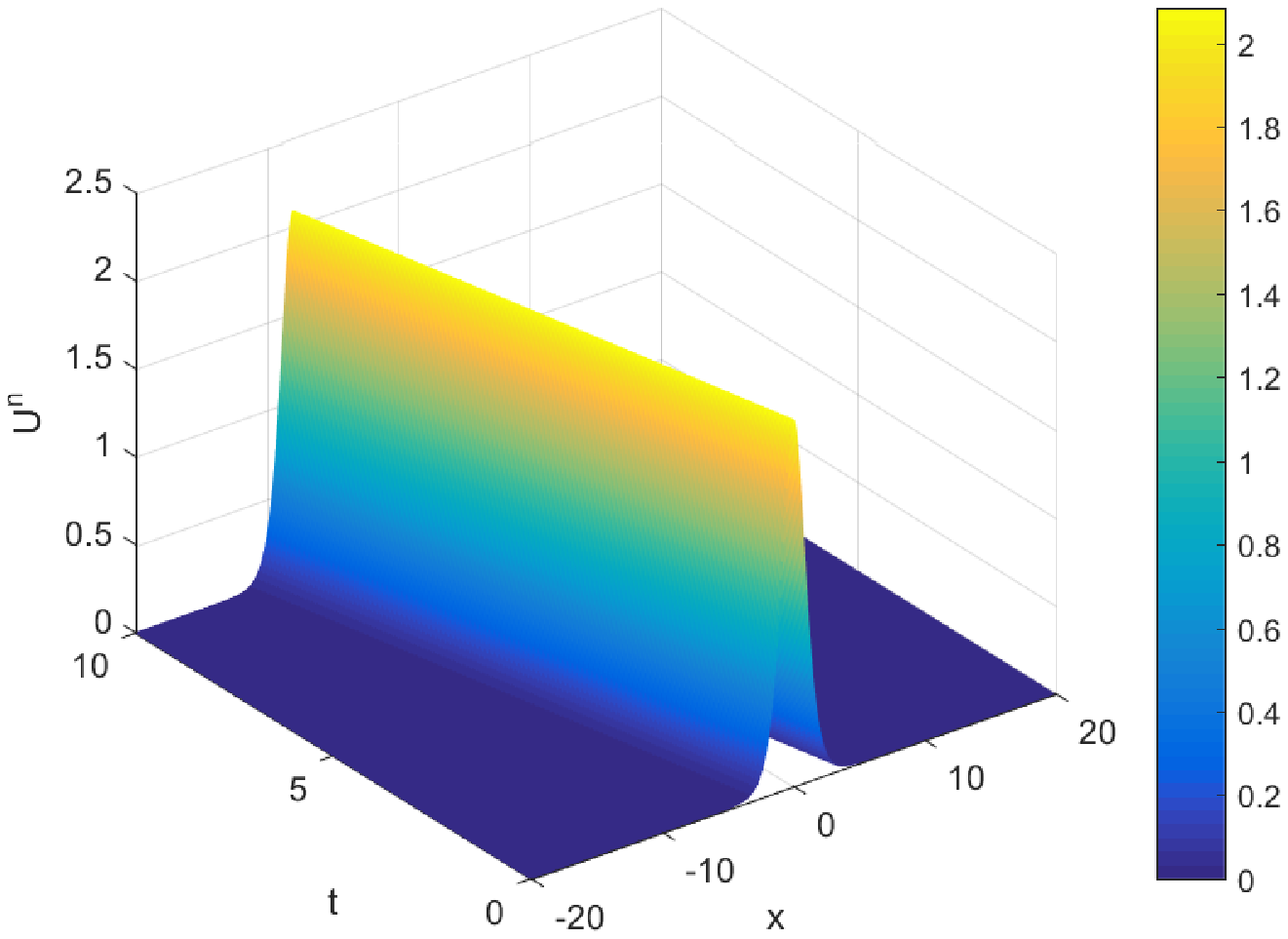}\\
{\footnotesize  \centerline {(e) $\alpha=2$}}
\end{minipage}
\begin{minipage}[t]{60mm}
\includegraphics[width=60mm]{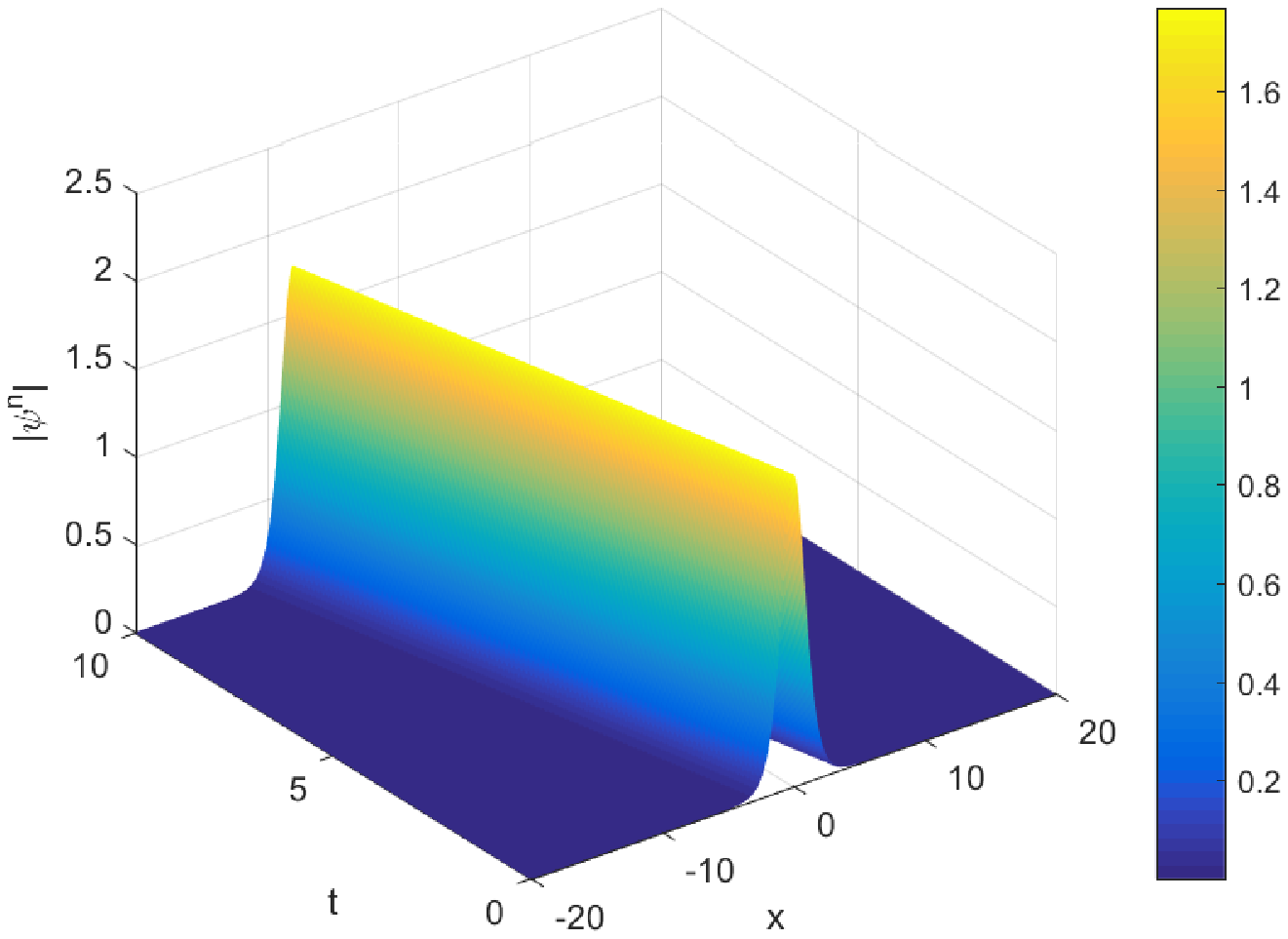}\\
{\footnotesize  \centerline {(f) $\alpha=2$}}
\end{minipage}
\caption{\small {\small{The evolution of the numerical solutions $U^n$  (left) and $|\varphi^n|$ (right) with different $\alpha$.} }}\label{fig522}
\end{figure}

\textbf{Example 4.2.} In the example, we study two-dimensional fractional KGS equation.
When $\alpha=\beta=2$, the exact solutions of this test problem are \begin{align*}
&\varphi(x,y,t)=A\exp\Big(\text{i}\big(\omega(x+y)-\theta t\big)\Big),\ \ u(x,y,t)=|\varphi(x,y,t)|,
\end{align*}
where $A=1, \omega=1, \theta=-0.5$. The region of this test problem is $[0,2\pi]\times[0,2\pi]$.
In order to measure the temporal error of our proposed methods, we choose $N_x=N_y=N=16$, which makes the spatial error small enough to ignore. In Fig. 7, we plot the errors versus $\tau$ at $T=1$. It shows that the FPAVF scheme is first-order accuracy in time, the FPAVF-C scheme and FPAVF-P scheme have the second-order temporal accuracy.
We depict the spatial errors of the three schemes with different $\alpha, \beta$ at $T=1$ in Fig. 8, which shows that the spatial error of the proposed schemes are very small and almost negligible, and the error
is dominated by the time discretization error. It confirms that, for sufficiently smooth
problems, the Fourier pseudo-spectral method is of arbitrary order in space.
In Fig. 9 and Fig. 10, we plot relative mass and energy errors of the three schemes. From the figures, we can see
that newly developed schemes can preserve the mass and energy conservation laws exactly.

\begin{figure}[H]
\centering
\begin{minipage}[t]{50mm}
\includegraphics[width=50mm]{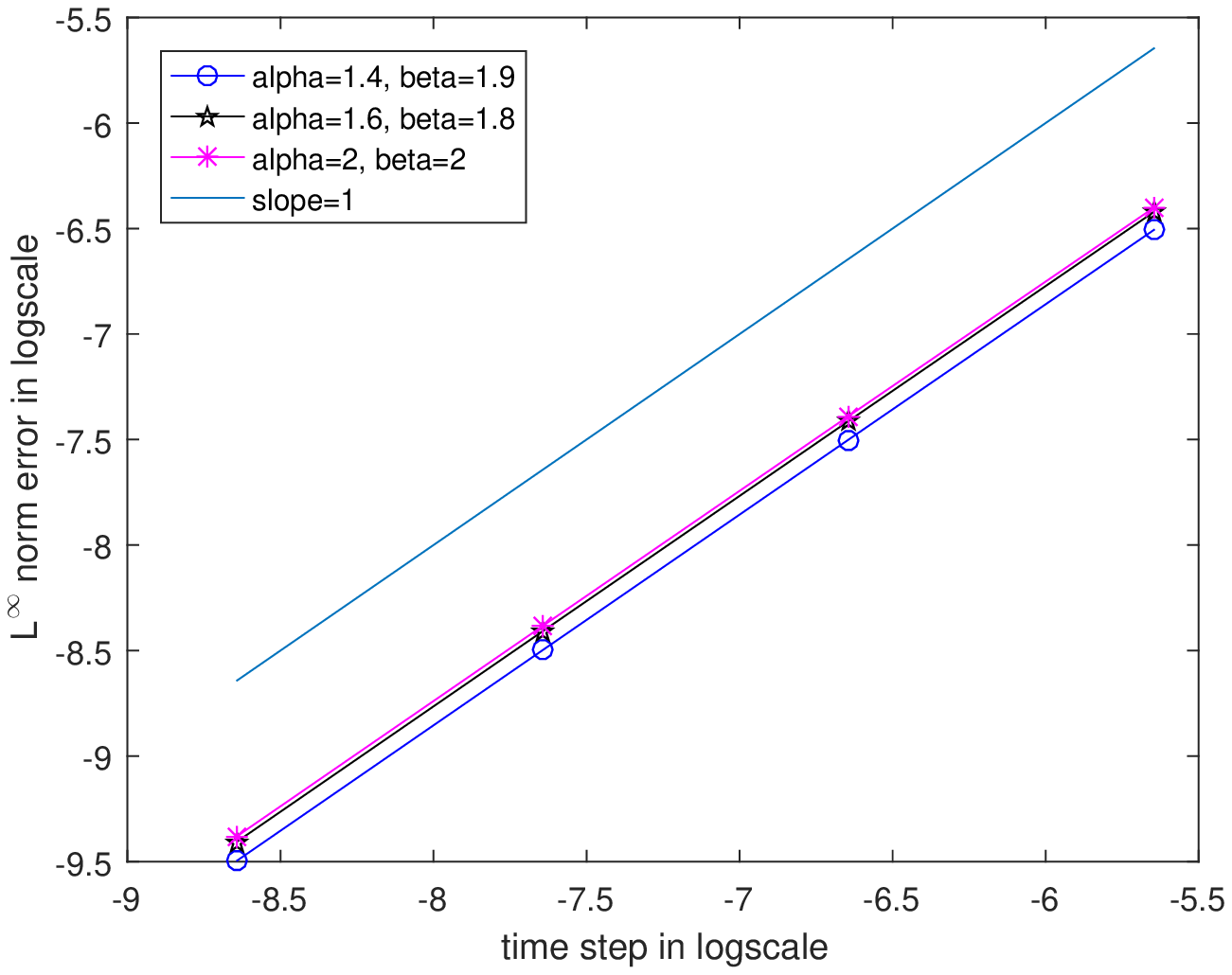}\\
{\footnotesize  \centerline {(a) FPAVF scheme}}
\end{minipage}
\begin{minipage}[t]{50mm}
\includegraphics[width=50mm]{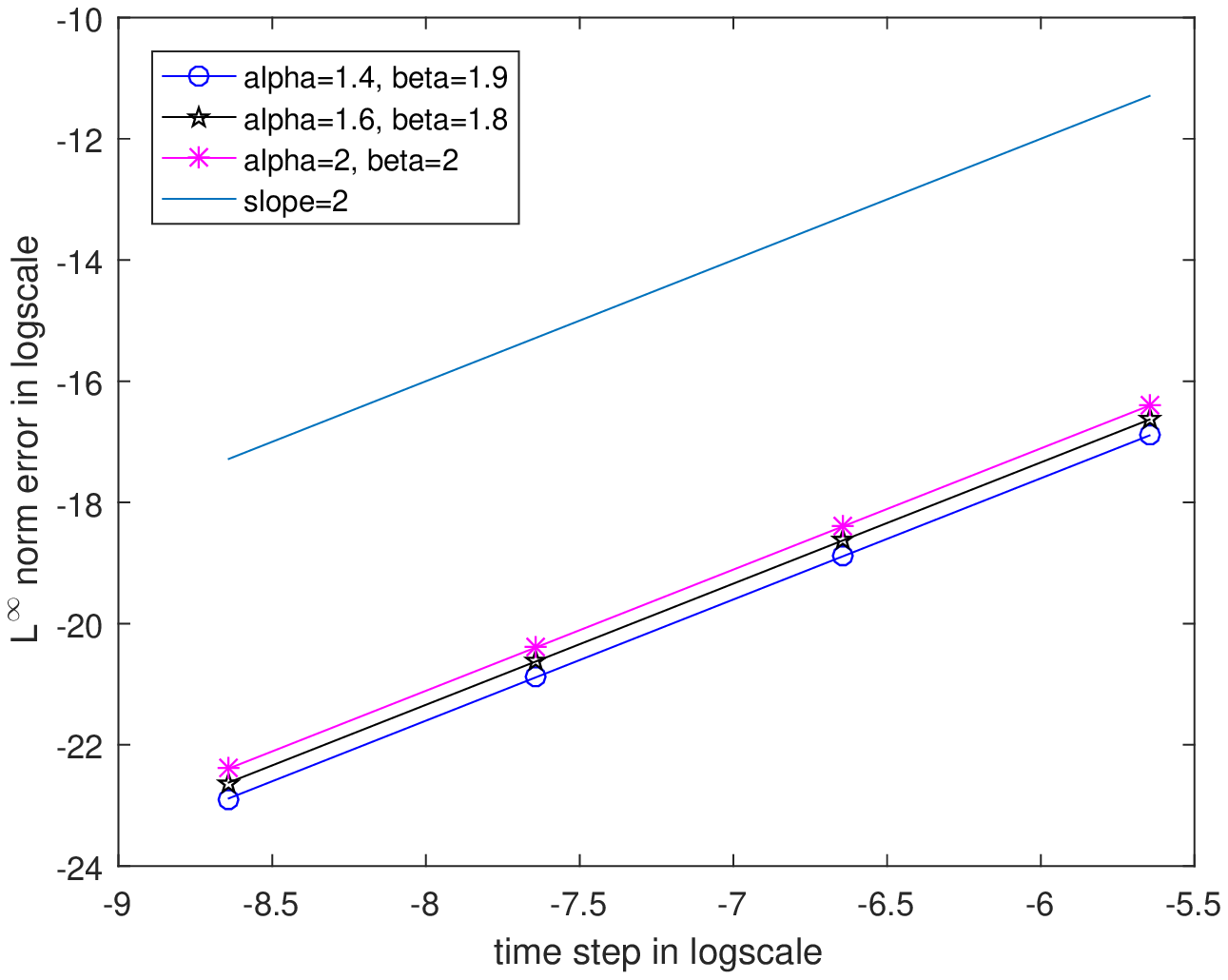}\\
{\footnotesize  \centerline {(b) FPAVF-C scheme}}
\end{minipage}
\begin{minipage}[t]{50mm}
\includegraphics[width=50mm]{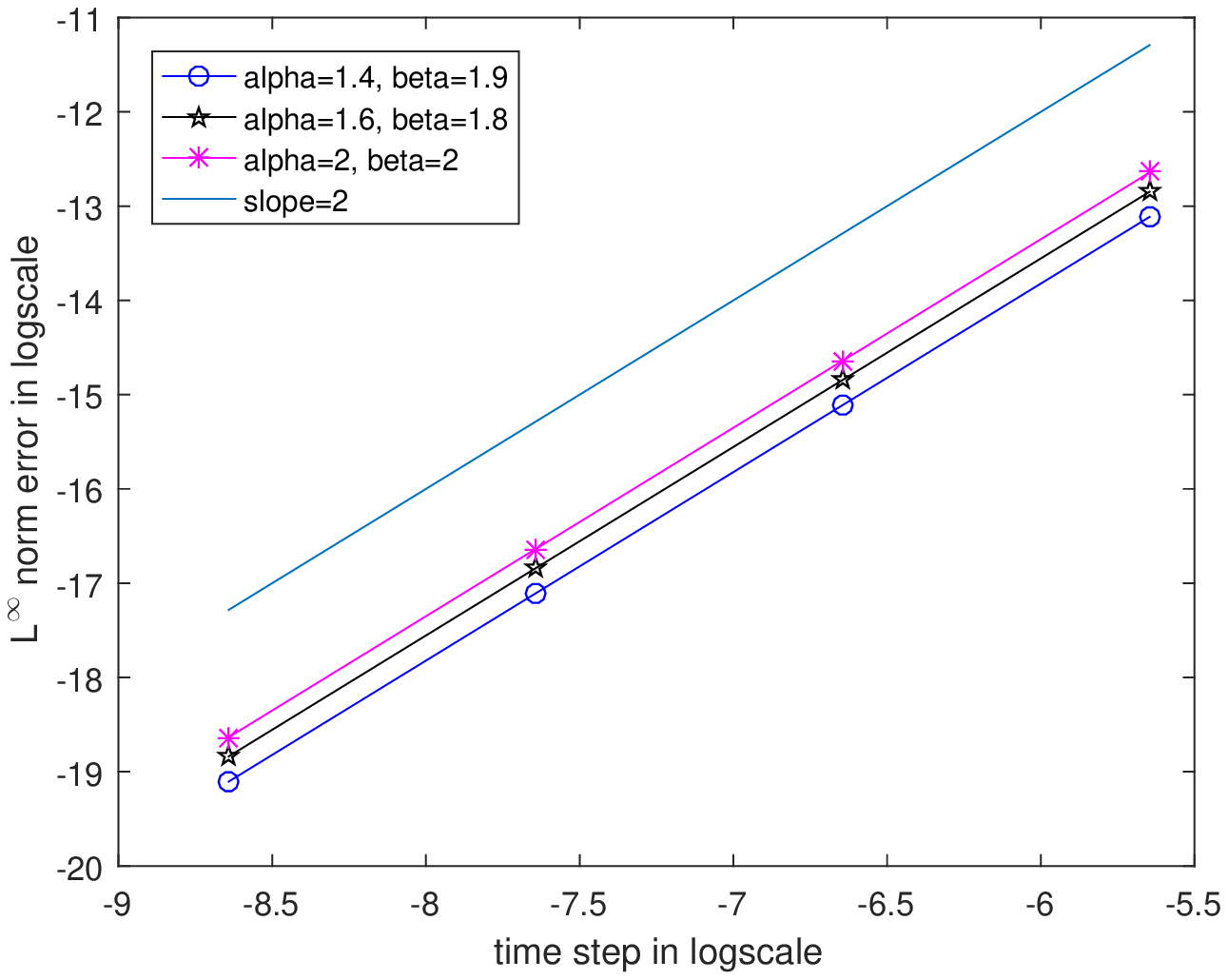}\\
{\footnotesize  \centerline {(c) FPAVF-P scheme}}
\end{minipage}
\caption{\small {Temporal accuracy of three schemes for different $\alpha$ and $\beta$ with $N=16$. }}\label{fig522}
\end{figure}

\begin{figure}[H]
\centering
\begin{minipage}[t]{50mm}
\includegraphics[width=50mm]{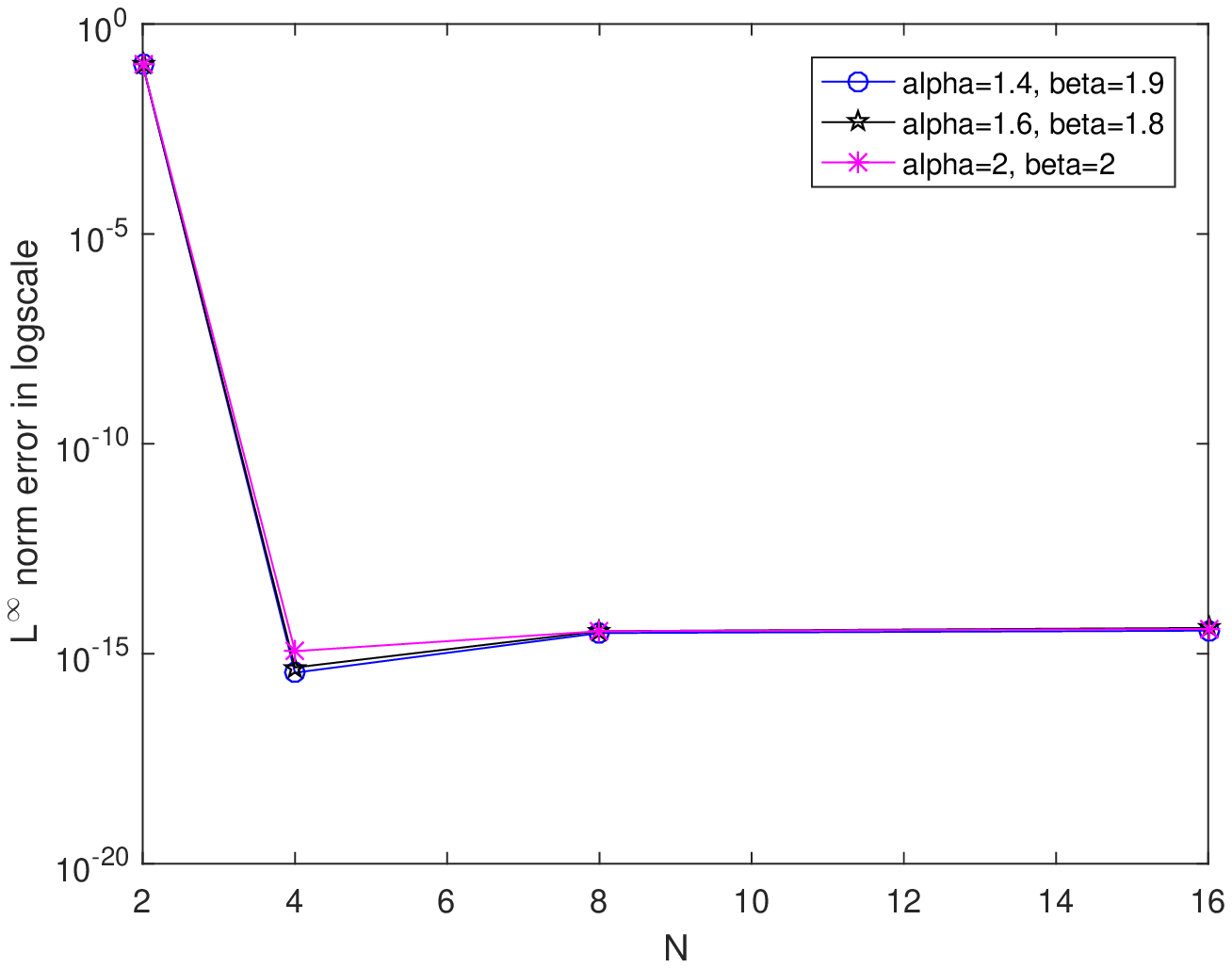}\\
{\footnotesize  \centerline {(a) FPAVF scheme}}
\end{minipage}
\begin{minipage}[t]{50mm}
\includegraphics[width=50mm]{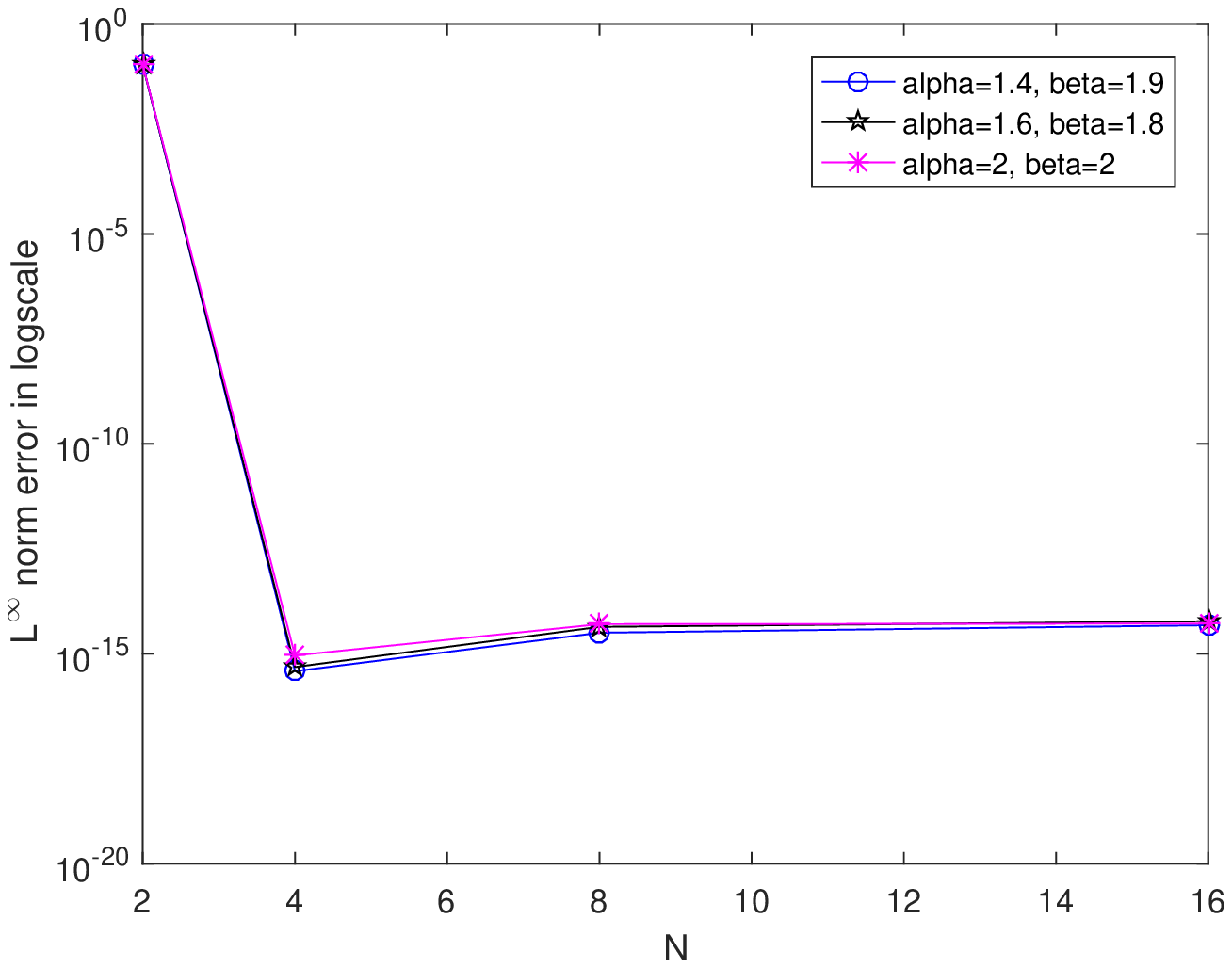}\\
{\footnotesize  \centerline {(b) FPAVF-C scheme}}
\end{minipage}
\begin{minipage}[t]{50mm}
\includegraphics[width=50mm]{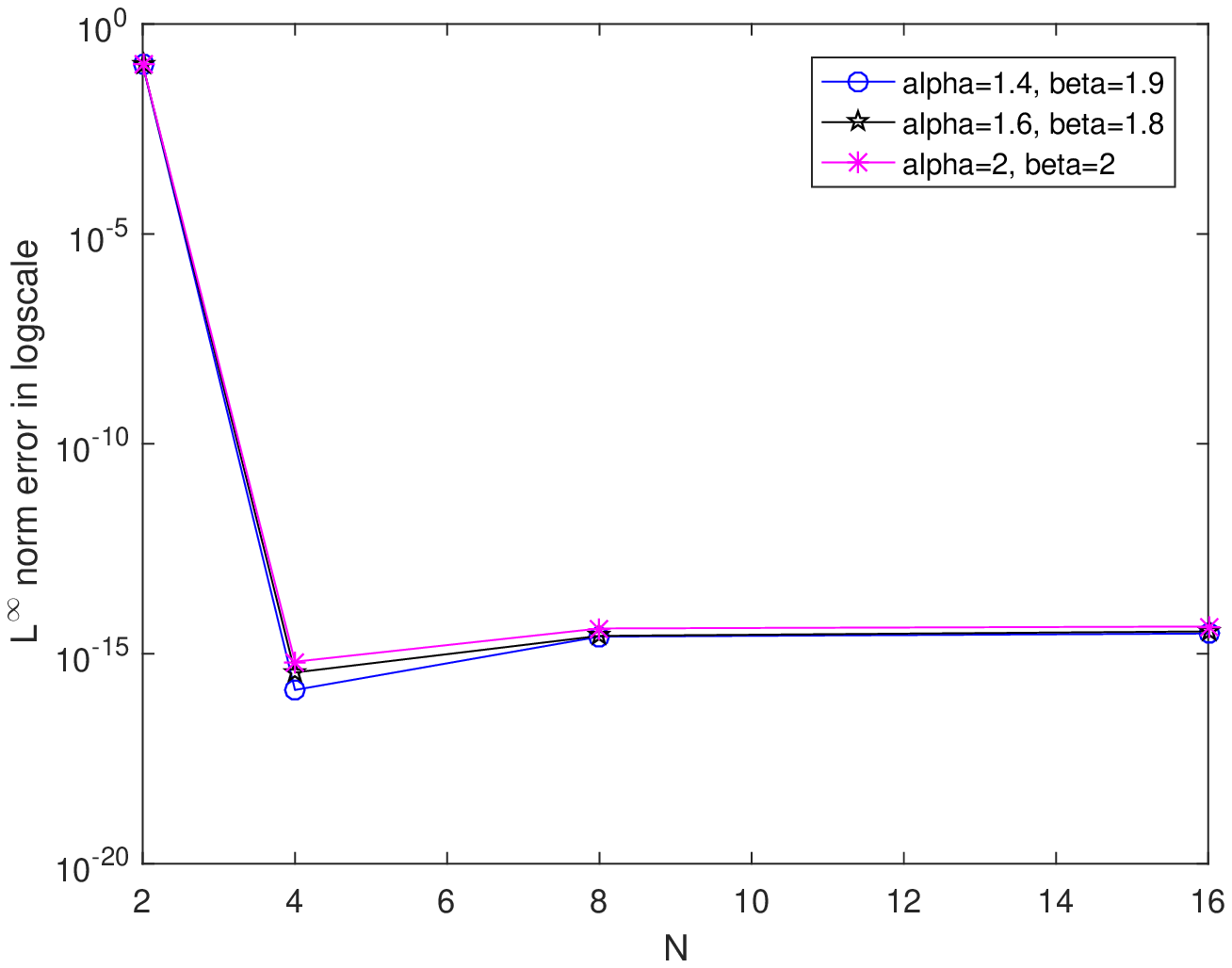}\\
{\footnotesize  \centerline {(c) FPAVF-P scheme}}
\end{minipage}
\caption{\small {Spatial accuracy of three schemes for different $\alpha$ and $\beta$ with $\tau=10^{-6}$.}}\label{fig522}
\end{figure}

\begin{figure}[H]
\centering
\begin{minipage}[t]{50mm}
\includegraphics[width=50mm]{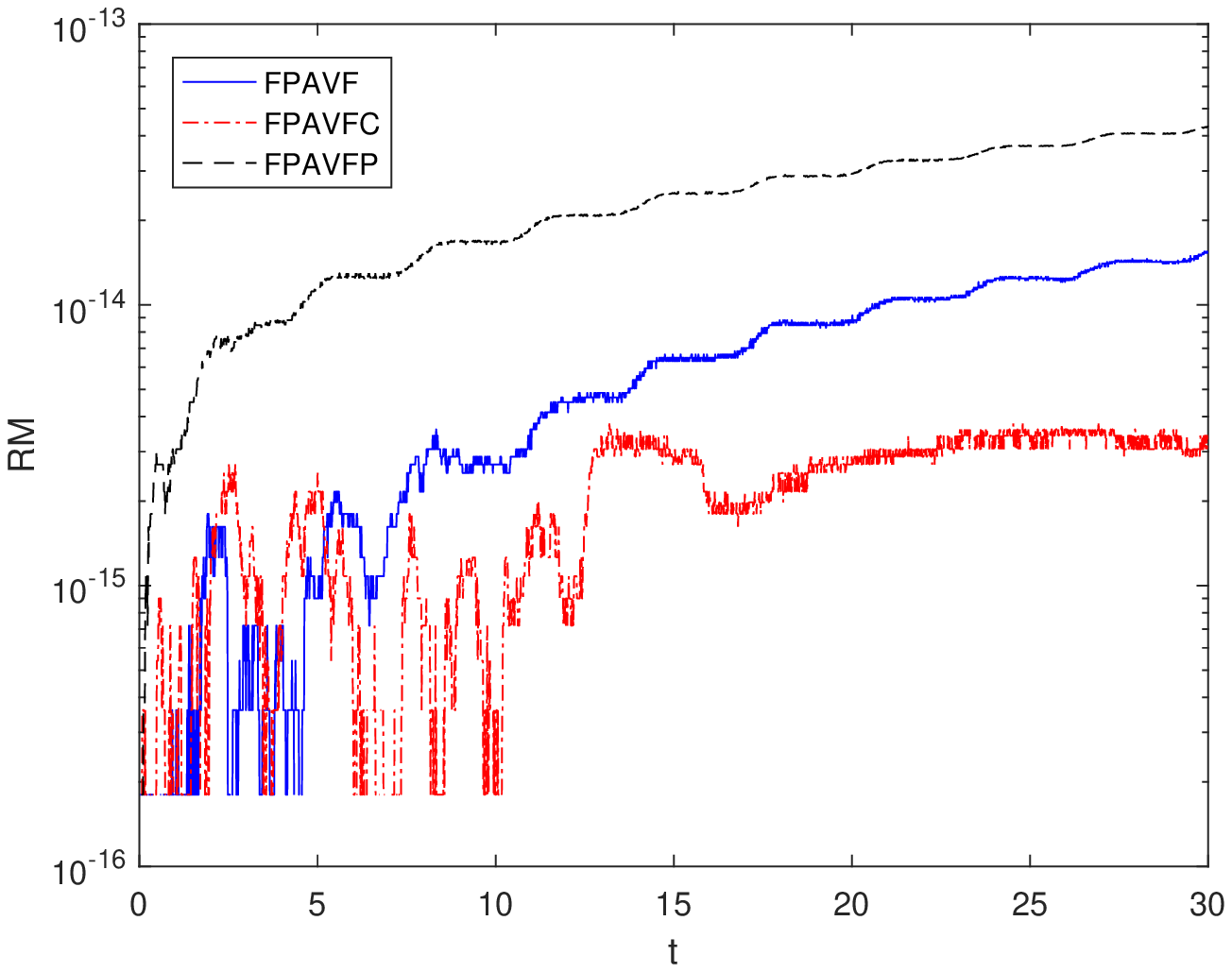}\\
{\footnotesize  \centerline {(a) $\alpha=1.4, \beta=1.9$}}
\end{minipage}
\begin{minipage}[t]{50mm}
\includegraphics[width=50mm]{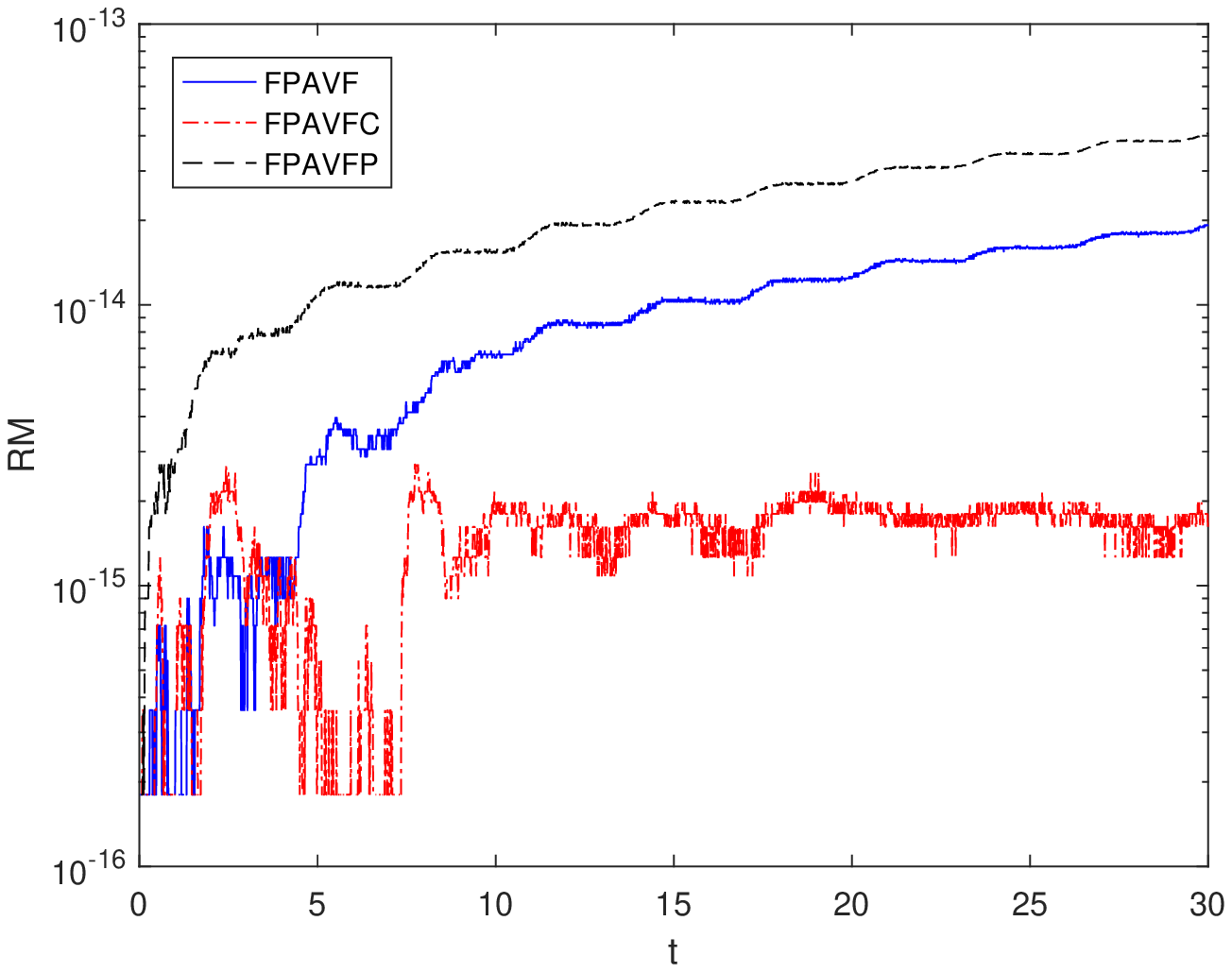}\\
{\footnotesize  \centerline {(b) $\alpha=1.6, \beta=1.8$}}
\end{minipage}
\begin{minipage}[t]{50mm}
\includegraphics[width=50mm]{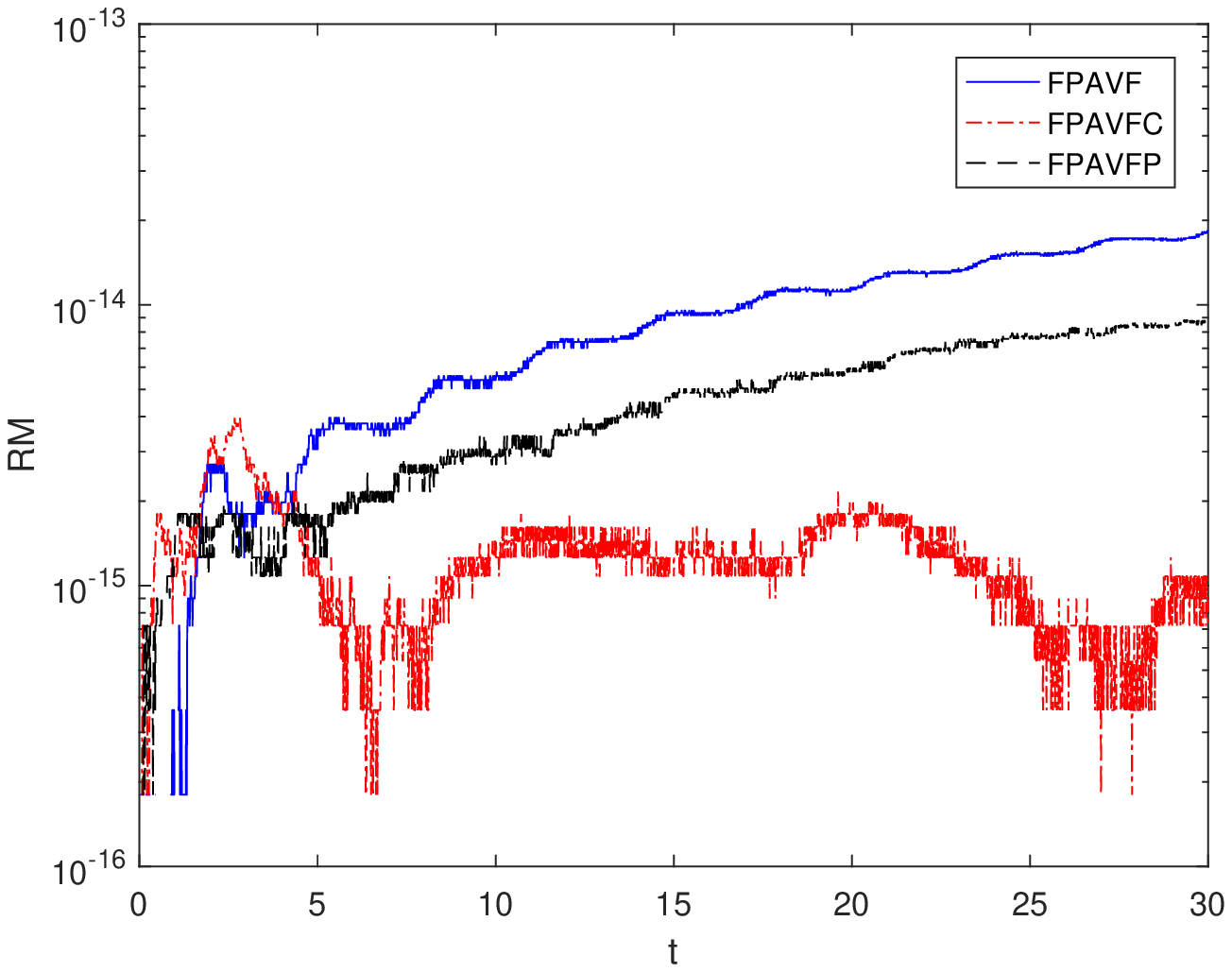}\\
{\footnotesize  \centerline {(c) $\alpha=2, \beta=2$}}
\end{minipage}
\caption{\small {The relative errors of mass with $N=16,~\tau=0.001$ for different $\alpha, \beta$. }}\label{fig522}
\end{figure}

\begin{figure}[H]
\centering
\begin{minipage}[t]{50mm}
\includegraphics[width=50mm]{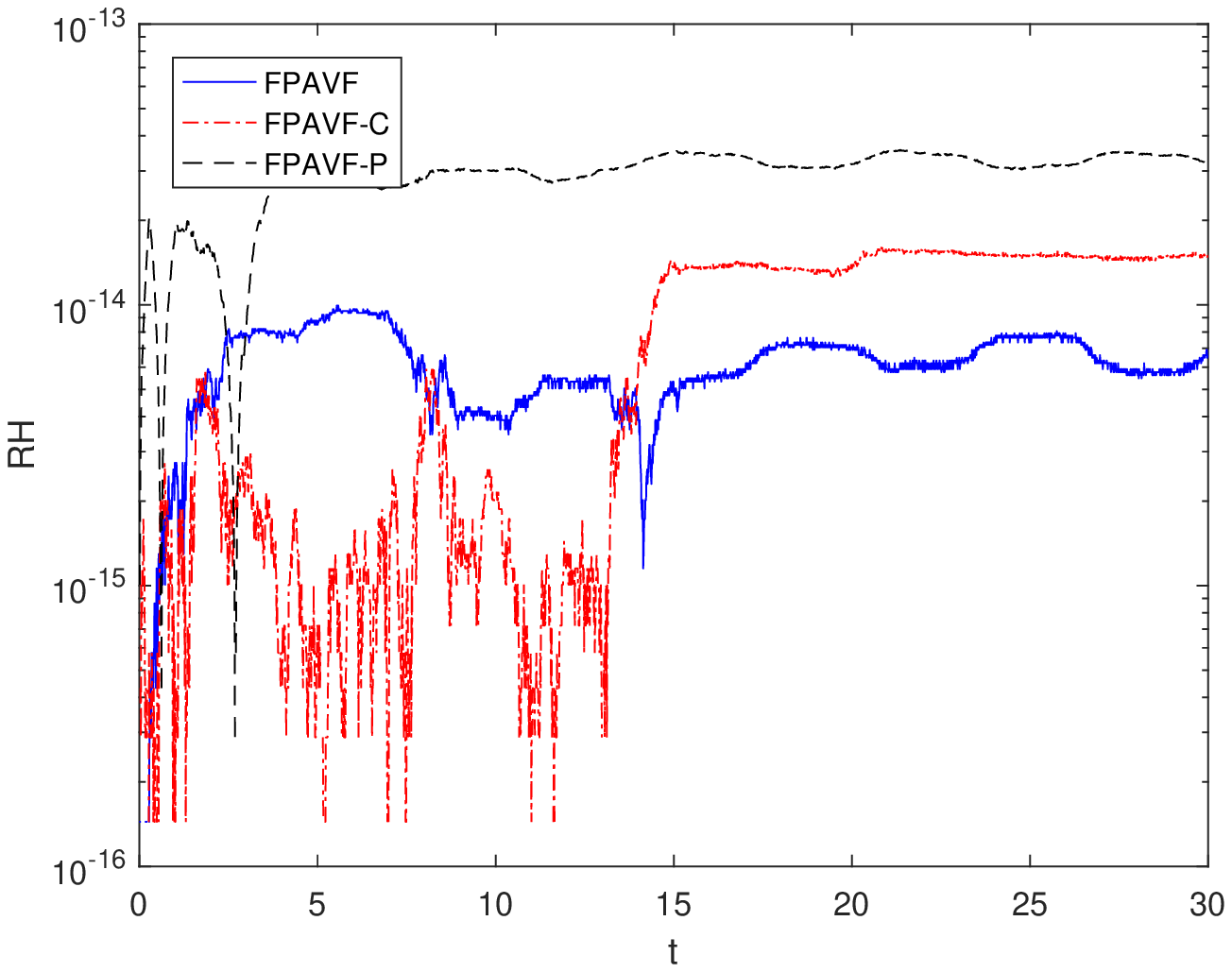}\\
{\footnotesize  \centerline {(a) $\alpha=1.4, \beta=1.9$}}
\end{minipage}
\begin{minipage}[t]{50mm}
\includegraphics[width=50mm]{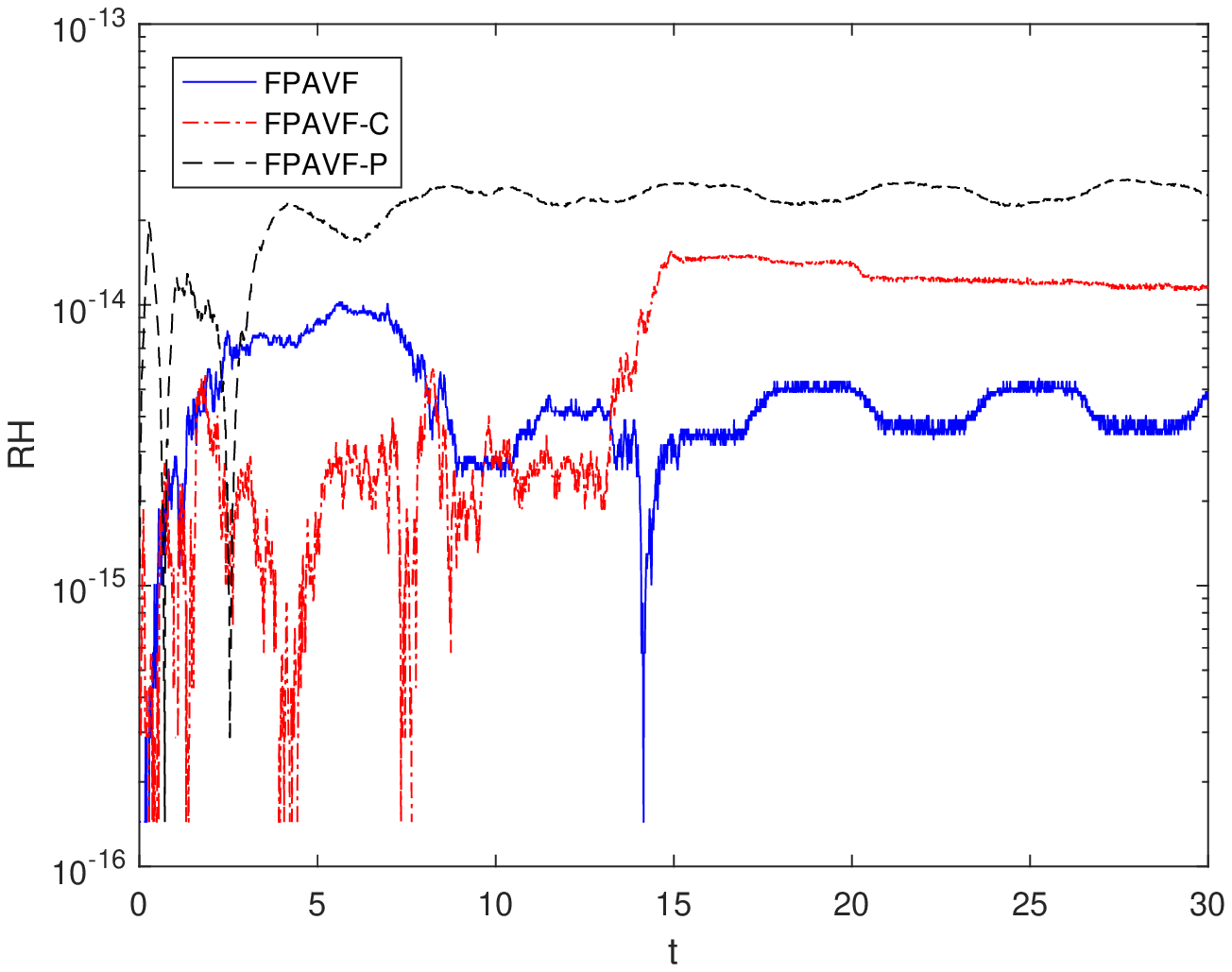}\\
{\footnotesize  \centerline {(b) $\alpha=1.6, \beta=1.8$}}
\end{minipage}
\begin{minipage}[t]{50mm}
\includegraphics[width=50mm]{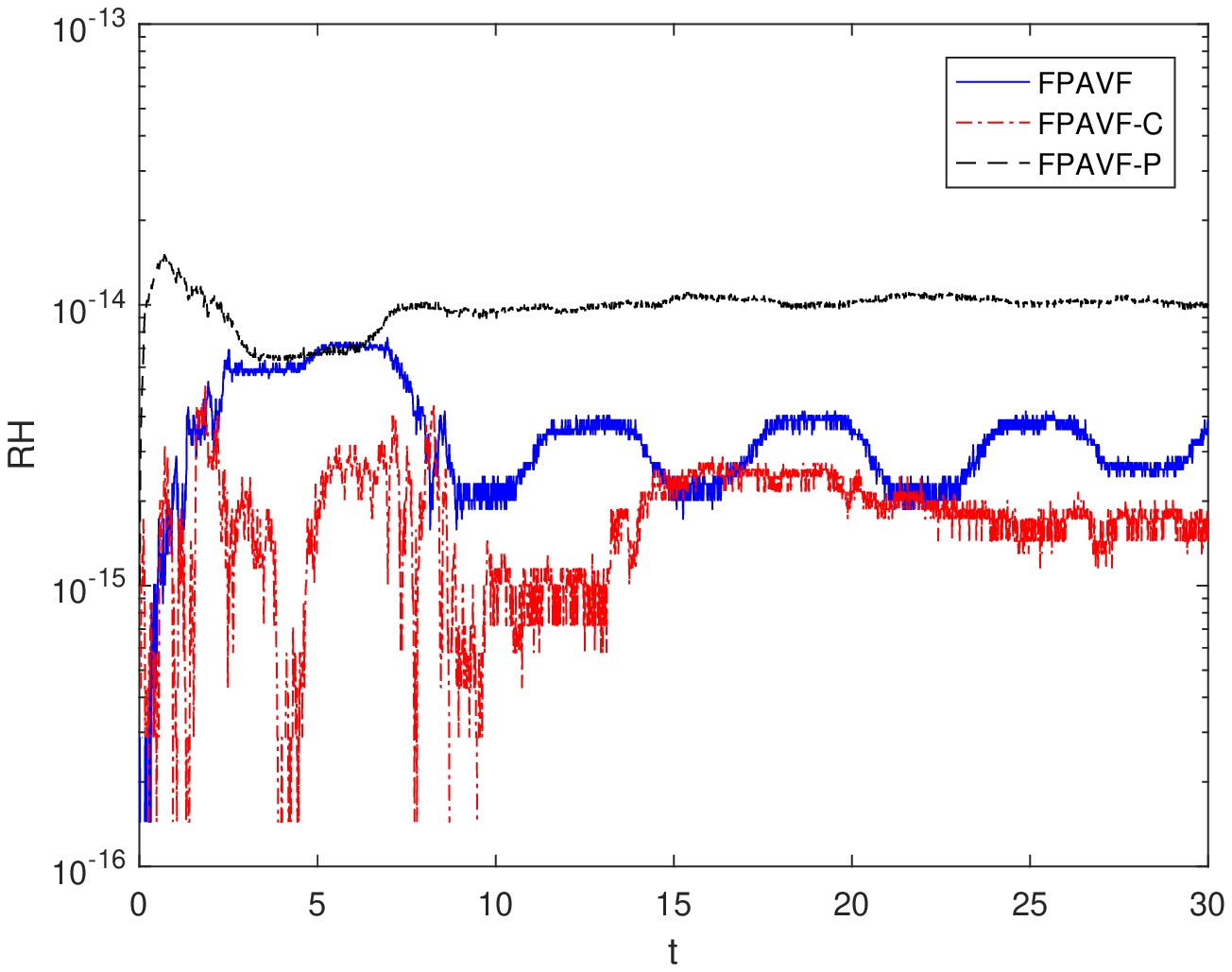}\\
{\footnotesize  \centerline {(c) $\alpha=2, \beta=2$}}
\end{minipage}
\caption{\small {The relative errors of energy with $N=16,~\tau=0.001$ for different $\alpha, \beta$.}}\label{fig522}
\end{figure}

\textbf{Example 4.3.}
Finally, we consider two-dimensional fractional KGS equations with the initial conditions
\begin{align*}
&\varphi(\bm x,0)=(1+\text{i})\exp(-|\bm x|^2),\ u(\bm x,0)=\text{sech}(|\bm x|^2),\ u_{t}(\bm x,0)=\sin(x+y)\text{sech}(-2|\bm x|^2),
\end{align*}
where $\bm x=(x,y)$. In this numerical example, we take $N_x=N_y=N=256$ and the temporal step size $\tau=0.001$, the computational domain
$\Omega=[-10,10]\times[-10,10]$. The relative errors in mass and energy are plotted in Fig. 11 and Fig. 12, which demonstrate that all schemes can conserve mass and energy conservation laws very well.
Fig. 13 shows the plots of $U^n$  and $|\varphi^n|$ with $\alpha=\beta=2$ at different times, which uniformly show a correct time evolutions compared to the results in Ref. \cite {p29}.

\begin{figure}[H]
\centering
\begin{minipage}[t]{60mm}
\includegraphics[width=60mm]{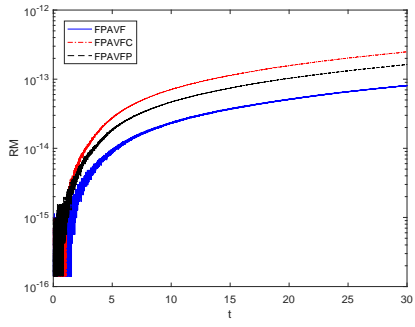}\\
{\footnotesize  \centerline {(a) $\alpha=1.5, \beta=1.8$}}
\end{minipage}
\begin{minipage}[t]{60mm}
\includegraphics[width=60mm]{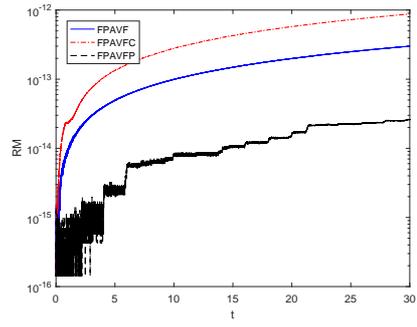}\\
{\footnotesize  \centerline {(b) $\alpha=2, \beta=2$}}
\end{minipage}
\caption{\small {The relative errors of mass for different $\alpha, \beta$. }}\label{fig522}
\end{figure}

\begin{figure}[H]
\centering
\begin{minipage}[t]{60mm}
\includegraphics[width=60mm]{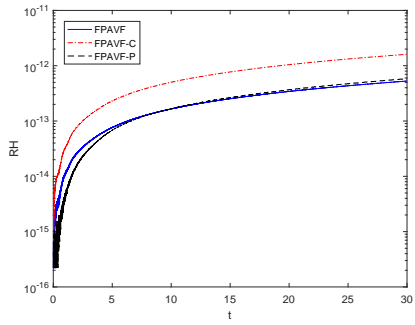}\\
{\footnotesize  \centerline {(a) $\alpha=1.5, \beta=1.8$}}
\end{minipage}
\begin{minipage}[t]{60mm}
\includegraphics[width=60mm]{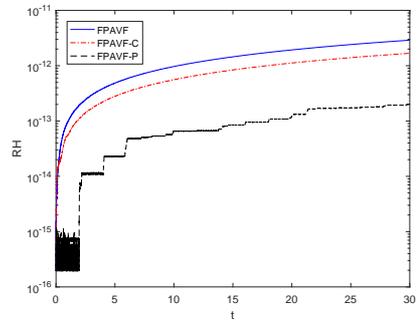}\\
{\footnotesize  \centerline {(b) $\alpha=2, \beta=2$}}
\end{minipage}
\caption{\small {The relative errors of energy  for different $\alpha, \beta$. }}\label{fig522}
\end{figure}

\begin{figure}[H]
\centering
\begin{minipage}[t]{60mm}
\includegraphics[width=60mm]{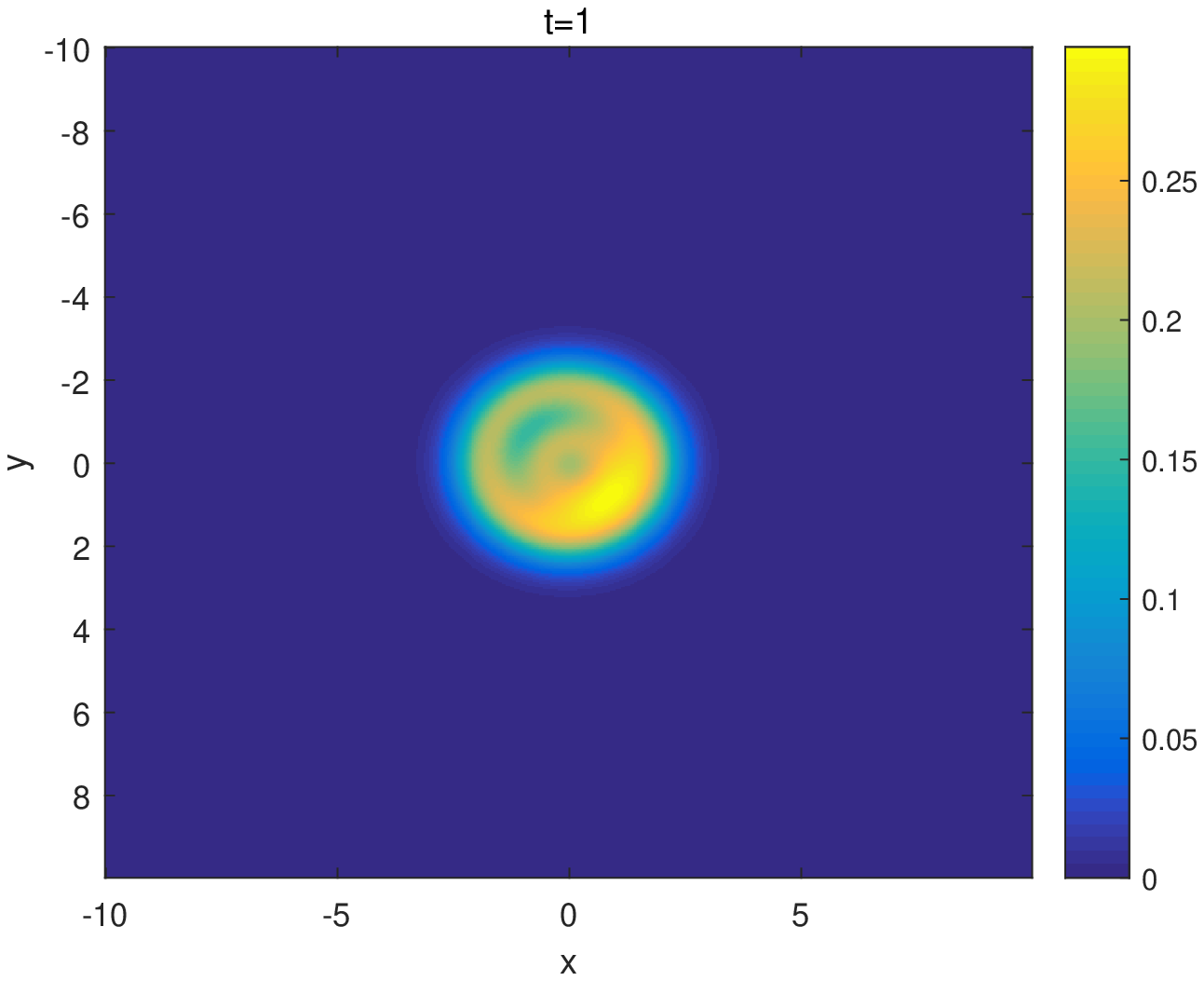}\\
\end{minipage}
\begin{minipage}[t]{60mm}
\includegraphics[width=60mm]{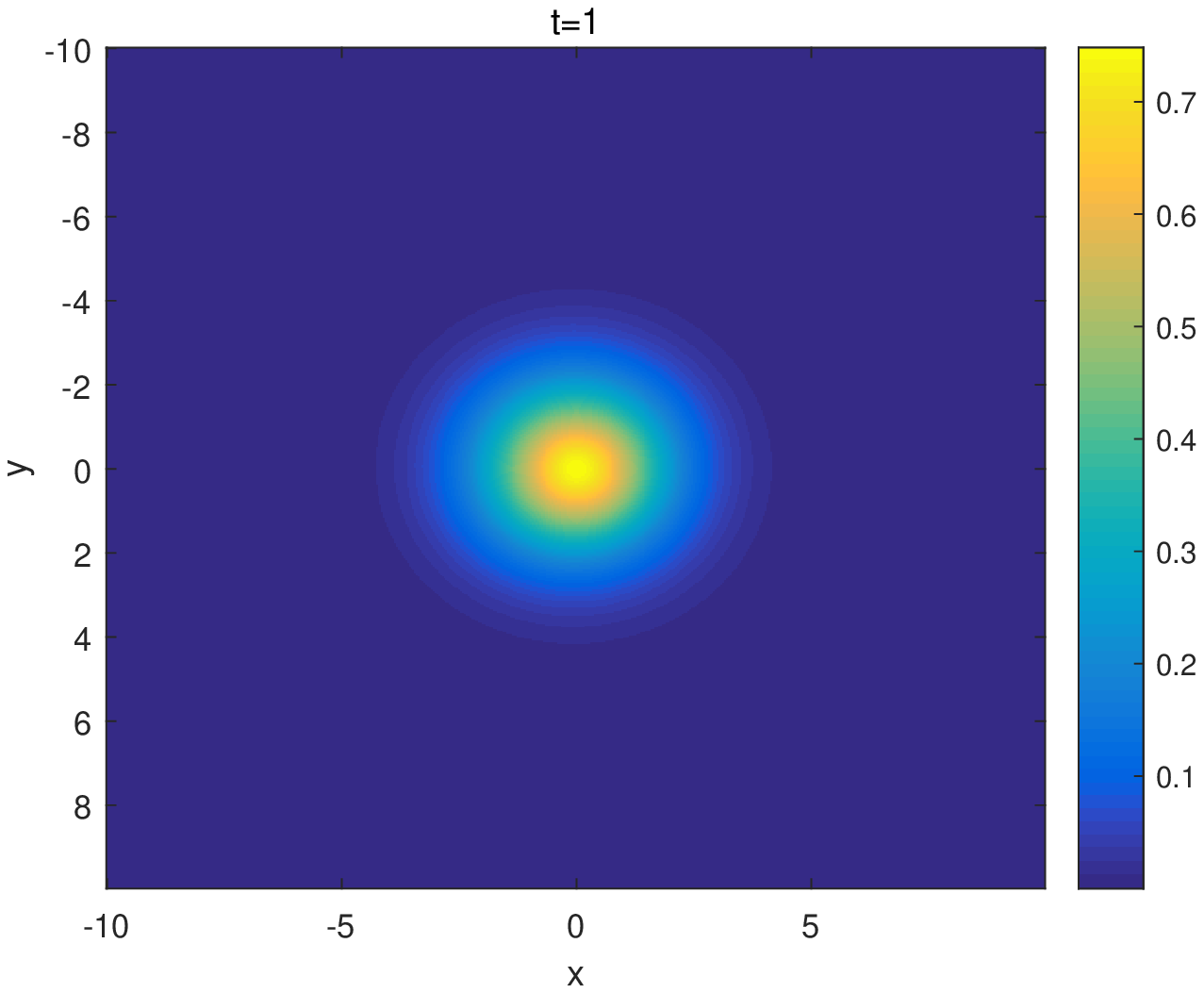}\\
\end{minipage}
\begin{minipage}[t]{60mm}
\includegraphics[width=60mm]{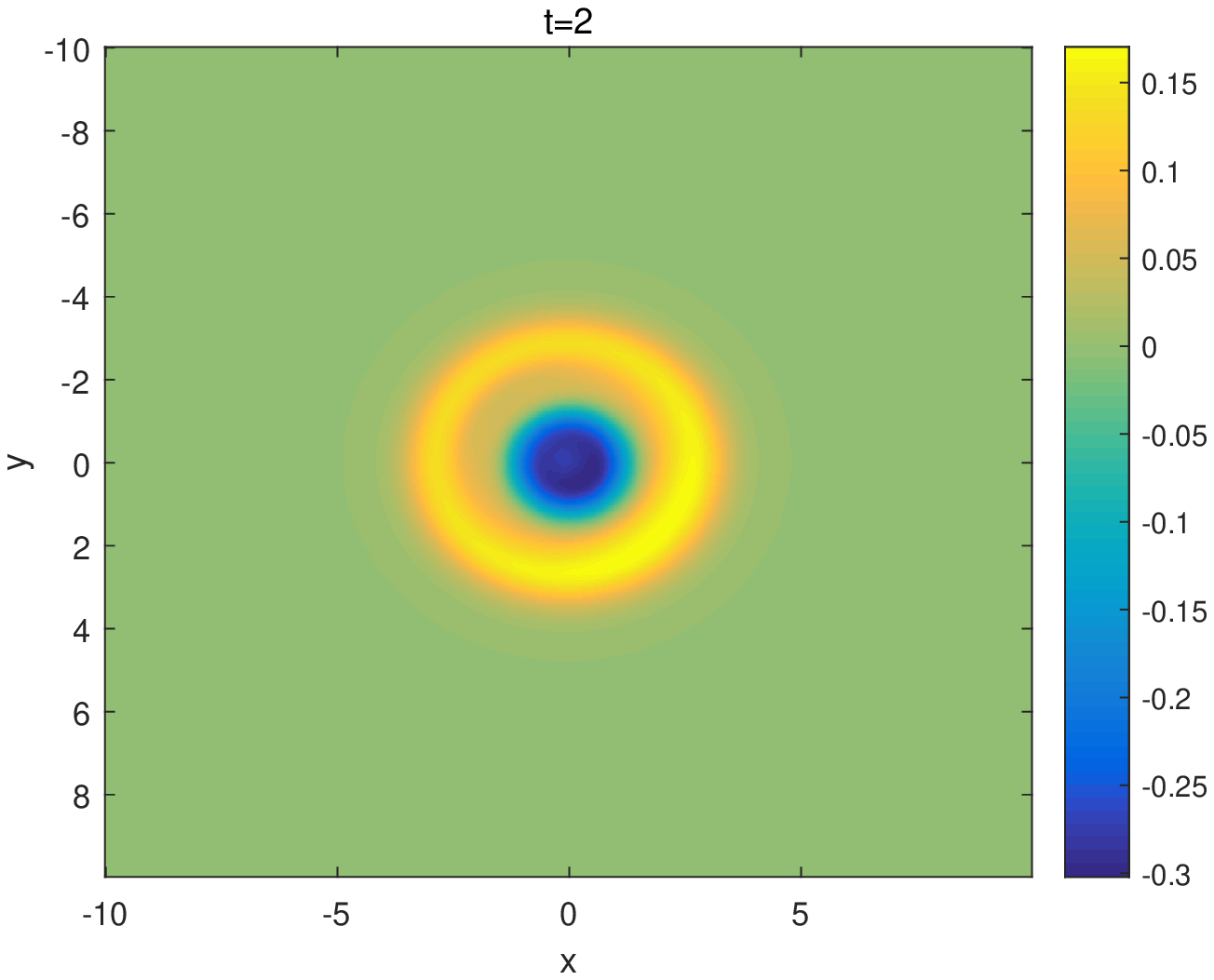}\\
\end{minipage}
\begin{minipage}[t]{60mm}
\includegraphics[width=60mm]{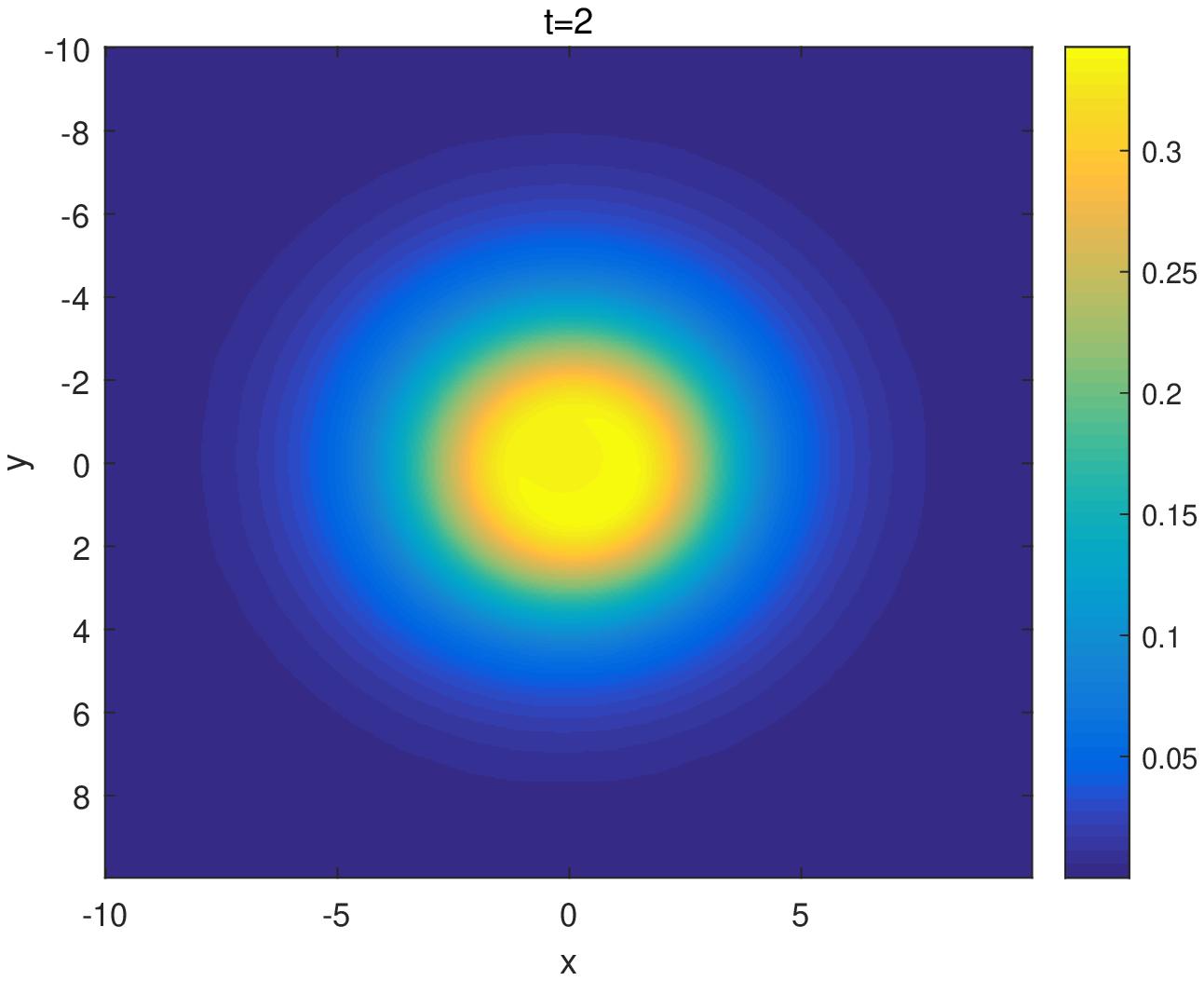}\\
\end{minipage}
\begin{minipage}[t]{60mm}
\includegraphics[width=60mm]{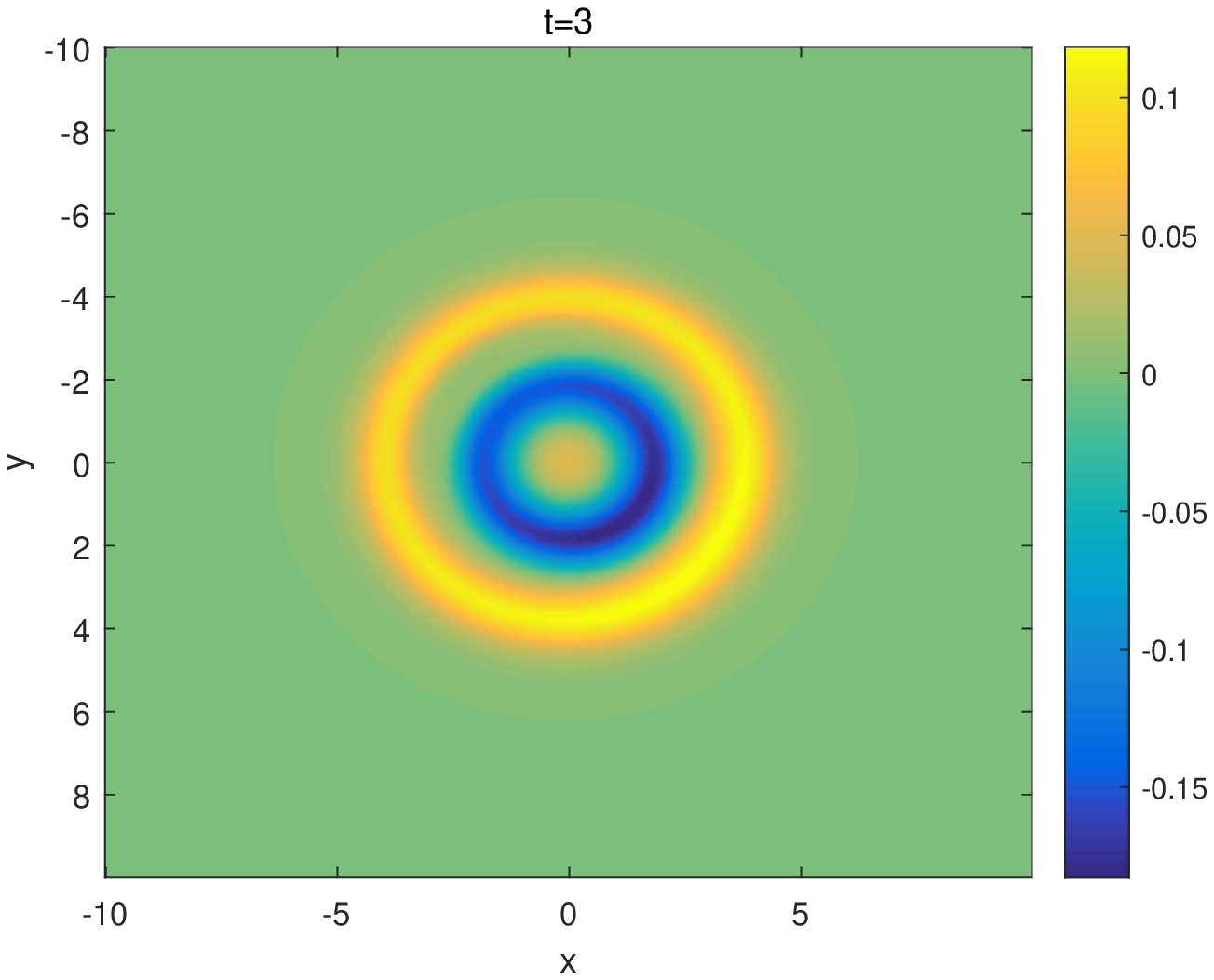}\\
\end{minipage}
\begin{minipage}[t]{60mm}
\includegraphics[width=60mm]{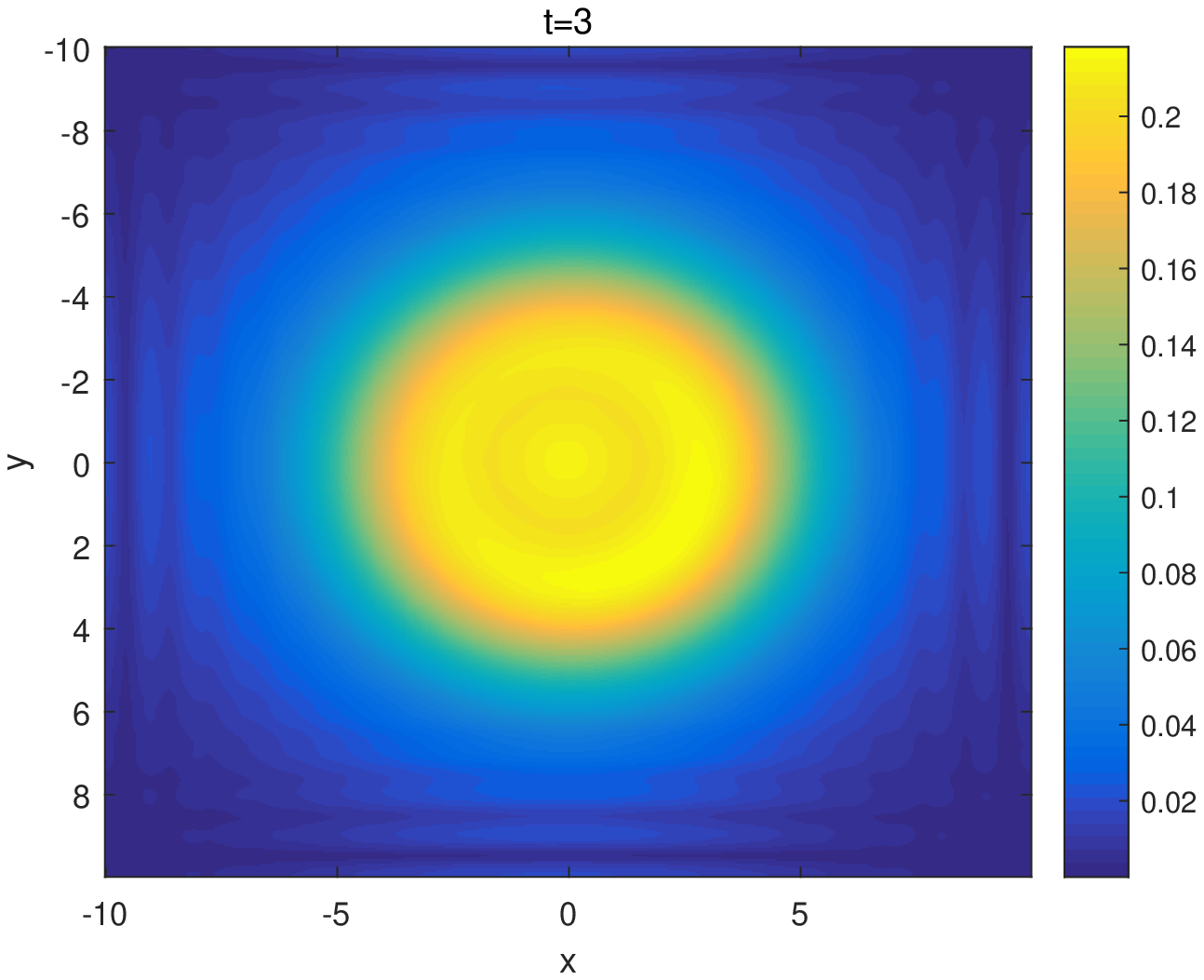}
\end{minipage}
\caption{\small {\small{Snapshot of $U^n$  (left) and $|\varphi^n|$ (right) by three conservative methods at $t=1, 2, 3$.} }}\label{fig522}
\end{figure}

\section{Conclusions and Remarks}

%In this paper, a new linear Fourier pseudo-spectral scheme, being built upon the scalar auxiliary variable approach, is developed to solve the two-dimensional fractional nonlinear Schr\"{o}dinger equation. The presented scheme is linear, accurate and can conserve energy very well. Theoretical analysis and numerical results show that the new scheme is not only much more efficient and easy to implement, it also preserve energy conservation law. Specifically, the scheme here can be carried over to other fractional equations, such as the fractional sine-Gordon equation, the fractional Klein-Gordon-Schr\"{o}dinger equation, etc.. Future research could continue to explore analysing errors of the schemes which are based on scalar auxiliary variable approach approach.
%is not only much more efficient and easy to implement, it also
%constructing some liner conservative schemes, based on the scalar auxiliary variable approach, for solving fractional equations, such as the fractional Klein-Gordon-Schr\"{o}dinger equation, the fractional nonlinear wave equations, etc.

In this paper, we reformulate the fractional Klein-Gordon-Schr\"{o}dinger equation as a canonical Hamiltonian system, and extend the  partitioned averaged vector field methods to construct a class of conservative schemes for the system based on its Hamiltonian formulation. Theoretical analysis and numerical results demonstrate that the proposed schemes are easy to implement and much more efficient, and can preserve both the discrete energy and mass conversation laws.

The Hamiltonian structure is important to analyse the conservative systems and further to construct numerical schemes for them.
The derivation of the Hamiltonian formulation is depended on the boundary conditions of the system. Up to now, only a few researchers consider the Hamiltonian formulation of fractional differential equations with periodic boundary conditions.
In our future work, we will investigate the properties of the fractional Laplacian and derive the Hamiltonian formulation of fractional differential equations with homogeneous boundary conditions, and develop structure-preserving algorithms based on the Hamiltonian formulation.
\section*{Acknowledgments}

This work is supported by the National Key Research and Development Project of China (Grant No. 2017YFC0601505, 2018YFC1504205),
the National Natural Science Foundation of China (Grant No. 11771213, 61872422, 11971242), the Postgraduate Research $\&$ Practice Innovation Program of Jiangsu Province (Grant Nos. KYCX19\_0776),
the Major Projects of Natural Sciences of University in Jiangsu Province of China (Grant No.18KJA110003), and the Priority Academic Program Development of Jiangsu Higher Education Institutions.

\end{document}